\documentclass[a4paper,10pt]{amsart}

\usepackage{amsmath,amsfonts,amsthm, amssymb,xypic, epsfig}
\usepackage[all]{xy}

\def\C{\mathbb{C}}
\def\Z{\mathbb{Z}}
\def\N{\mathbb{N}}

\def\qed{$\hfill \diamondsuit$}
\def\U{\mathbf{U}}
\def\H{\mathbf{H}}
\def\h{\mathfrak{h}}
\def\g{\mathfrak{g}}
\def\n{\mathfrak{n}}

\def\O{{\mathcal{O}}}

\def\Dtn{\Delta^{(1)}_n}
\def\b{\mathbf{b}}
\def\f{\mathbf{f}}
\def\m{\mathbf{m}}
\def\d{\mathbf{d}}
\def\AH{\widehat{\mathbf{H}}}
\def\Si{\widehat{\mathfrak{S}}}
\def\F{\mathbb{F}_q}

\newtheorem{theo}{Theorem}
\newtheorem{prop}{Proposition}

\newtheorem{cor}{Corollary}

\newtheorem*{conj}{Conjecture}
\numberwithin{equation}{section}
\numberwithin{lem}{section}
\numberwithin{prop}{section}
\numberwithin{cor}{section}
\numberwithin{theo}{section}
\numberwithin{claim}{section}

\begin{document}

\title{ Quivers of type $A$, flag varieties and representation theory}
\author{Olivier Schiffmann}
\date{}

\maketitle
\noindent

\paragraph{\textbf{Introduction.}}In this survey, we describe and 
relate various occurences of quivers
of type $A$ (both
finite and affine) and their canonical bases in combinatorics, in algebraic geometry and in 
representation theory.
The ubiquity of these quivers makes them especially important to 
study : they are
pervasive in very classical topics (such as the theory of symmetric 
functions) as well
as in some of the most recent and exciting areas of representation theory
(such as representation theory of quantum affine algebras, affine 
Hecke algebras or
Cherednik algebras).

\vspace{.2in}

\paragraph{}There is a vast litterature on the subject and we have 
been forced to
make a choice in the selection of the results presented here.
We believe that one of the main reasons for the omnipresence of 
quivers of type $A$
is the fact that they are related (in more than one way) to classical and
fundamental objects in geometric representation theory of type $A$ such as
(partial) flag varieties and nilpotent orbits. This is the point we tried to
emphasize in this survey. For that reason,
many interesting and important results are only alluded to or sketched here and
we apologize to all those whose work we did not mention through a lack of 
space and competence.
Also, we do not give complete proofs of all results presented here and refer to
the original papers whenever possible. Rather, we have tried to 
convey the fundamental
ideas in a nontechnical way as much as possible.

\vspace{.2in}

\paragraph{}This survey is organized as follows. In Section~1 we describe some basic geometric
properties of quivers of type $A$ and their representation spaces. Section~2 is a brief
reminder on Ringel-Hall algebras and their canonical bases,
and in Sections~3 and 4 we describe the Ringel-Hall algebras of the
categories of representations of quivers of type $A_1^{(1)}$ and type $A_n, A_{n}^{(1)}$
respectively. This makes a link with combinatorics (symmetric functions, \cite{Mac}) and with 
quantum
group theory (\cite{R}, \cite{S1}). 
In Section~5 we recall some classical geometric constructions in type $A$
such as (partial, affine) flag varieties, nilpotent orbits, Hecke algebras and
various convolution
algebras. In Section~6, we give a first link between quivers of type $A$ and (partial, affine)
Schubert varieties; one important consequence is that the coefficients of the canonical basis
of the Hall algebras of quivers of type $A$ are given by Kazhdan-Lusztig polynomials 
(\cite{Lus1},
\cite{Zel}).
Section~7 gives the geometric realization of Schur-Weyl duality and its affine and
quantum cousins. This gives a conceptual explanation for some of the results in Section~4.
We totally change perspective in Section~8 and, following \cite{LLT}, \cite{Ar} and 
\cite{VV}, 
we present a surprising and beautiful application
of the theory of
canonical bases and quivers of type $A$ to character forumlae for quantum affine algebras
and affine Hecke algebras (at roots of unity). Finally, we gather a few more interesting
occurences of quivers of type $A$ in relation to geometry (Nakajima quiver varieties,
parabolic structures on curves and surfaces) in the last Section.

\section{Quivers of type $A$}
\paragraph{\textbf{1.1.}} We will denote by $\Delta_n$ the 
(equioriented) quiver of type $A_n$, with vertex set
$I_n=\{1, \ldots, n\}$ and with edge set $\{(i,i+1)\;|\;i=1, \ldots, 
n-1\}$. 
Similarly,
we let
$\Dtn$ stand for the (equioriented) affine quiver of type 
$A_{n-1}^{(1)}$, with vertex set $\Z/n\Z$ and edge
set $\{(i,i+1)\;|i \in \Z/n\Z\}$ (note that this also makes sense for $n=1$). 
Let $k$ be any field. Recall that a 
\textit{representation of $\Delta_n$} (resp.
\textit{of $\Dtn$}) \textit{over $k$}
consists of a pair $(\mathbf{V},\underline{\phi})$ where 
$\mathbf{V}=\bigoplus_{i }  V_i$ is an $I_n$-graded
(resp. $\Z/n\Z$-graded) $k$-vector space and 
$\underline{\phi}=(\phi_{i}) \in \bigoplus_i
\mathrm{Hom}(V_i, V_{i+1})$.
We let $Mod_k(\Delta_n)$ be the abelian category of representations 
of $\Delta_n$ over $k$.  In the case of $\Dtn$
we will impose an additional \textit{nilpotency} condition; namely, we 
denote by $Mod_k(\Dtn)$ the abelian category
of representations $(\mathbf{V},\underline{\phi})$ of $\Dtn$ such 
that $\phi_N \cdots \phi_2 \phi_1=0$ for $N \gg 0$.
This restriction may seem artificial from the strict point of view of 
quivers, but it will appear naturally later in the
context of flag varieties.

\vspace{.2in}

\paragraph{\textbf{1.2.}} Let $\Delta$ be either $\Delta_n$ or 
$\Dtn$. The simple objects in $Mod_k(\Delta)$ are
naturally in bijection with the vertices of $\Delta$, and we will 
denote by $S_i$ the simple module with
$\mathbf{V}=V_i= k$.  Indecomposable objects of 
$Mod_k(\Delta_n)$ are parametrized by pairs
$(i;l)$ with $i \in I_n$ and $1 \leq l \leq n+1-i$ : to $(i;l)$ 
corresponds the unique indecomposable object
$S_{i,l}$ with top $S_i$ and length $l$. Similarly, indecomposable 
objects of $Mod_k(\Dtn)$ are parametrized by
pairs $(i;l)$ with $i \in \Z/n\Z$ and $l \geq 1$, and we use the same 
notation $S_{i,l}$ for the unique indecomposable
with top $S_i$ and length $l$. In particular, in both cases, the set 
of isomorphism classes of objects in
$Mod_k(\Delta)$ is in bijection with a set of 
\textit{n-multipartitions} (or \textit{n-multisegments})
$$\mathbf{m}=(\underline{\lambda}^{(i)})_i=\sum_i \sum_t (i;\lambda^{(i)}_t)$$
where $\lambda^{(i)}_t \leq n+1-i$ if $\Delta=\Delta_n$.

\vspace{.2in}

\paragraph{\textbf{1.3.}}  Let $\Delta$ again stand for either 
$\Delta_n$ or $\Dtn$.
Let us fix a dimension vector
$\mathbf{d}=(d_i)_i \in \N^{I}$ with
$I=I_n$ or $I=\Z/n\Z$ according to $\Delta$. The group 
$G_{\mathbf{d}}=\prod_i GL_k(d_i)$ acts naturally by
conjugation on the space
$$E(\mathbf{d})=\bigoplus_{i=1}^{n-1} \mathrm{Hom}(k^{d_i},k^{d_{i+1}})$$
if $\Delta=\Delta_n$ and on the space
$$E(\mathbf{d})=E_{nilp}(\mathbf{d})=\{(\phi_i) \in \bigoplus_{i \in \Z/n\Z} \mathrm{Hom} 
(k^{d_i}, k^{d_{i+1}})\:|\phi_N \cdots \phi_2 \phi_1=0
\;\mathrm{for}\; N \gg 0\}$$
if $\Delta=\Dtn$. Moreover, the set of $G_{\mathbf{d}}$-orbits is (by definition) in 
bijection with the set of isomorphism classes of objects
in $Mod_k(\Delta)$ of dimension $\mathbf{d}$. We will denote by 
$\mathcal{O}_{\mathbf{m}}$ the
$G_{\mathbf{d}}$-orbit (in the appropriate space) corresponding to a 
multipartition $\mathbf{m}$.\\
\hbox to1em{\hfill}Since there are only finitely many multipartitions 
of a given dimension $\mathbf{d}$ the
group $G_{\mathbf{d}}$ acts on $E(\mathbf{d})$ with only finitely many orbits.
It follows that there is a (unique) dense orbit (the ``generic 
representation of dimension $\mathbf{d}$''); it is also
easy to see that there is a unique closed orbit (the semisimple 
representation of dimension $\mathbf{d}$, in which all maps $\phi_i$ are zero).
Observe that
when $\Delta=\Delta_1^{(1)}$ we recover the classical situation of 
the group $GL(d)$ acting on the nilpotent cone
$\mathcal{N}=\{x \in gl_d\;|x^N=0\;\mathrm{for}\; N \gg 0\}$.

\vspace{.2in}

\paragraph{\textbf{1.4.}} We use the same notation as in 
\textbf{1.3.} Let $\mathbf{m}$ be a multisegment
of dimension $\mathbf{d}$. The Zariski closure 
$\overline{\mathcal{O}_{\mathbf{m}}}$ is a union of
$G_{\mathbf{d}}$-orbits. This induces a partial order $\leq$ on the set of 
orbits, which can be explicitly described.
More important for us will be the geometric properties of 
the singular variety
$\overline{\mathcal{O}_{\mathbf{m}}}$. The first result in this 
direction is the following.
\begin{theo}[\cite{ADFK},\cite{LM}] The variety $\overline{\mathcal{O}_{\mathbf{m}}}$ is normal.
\end{theo}
This theorem has been generalized to quivers of type $A$ and $D$ with 
arbitrary orientations by
Bobinski-Zwara (\cite{BZ1}, \cite{BZ2}), and is conjectured to hold for arbitrary 
quivers (in the case of the nilpotent cone this is a well-known result of
Kostant). In the same spirit,  Abeasis, Del Fra and Kraft constructed in \cite{ADFK}
an explicit desingularization of 
$\overline{\mathcal{O}_{\mathbf{m}}}$ as follows :
\paragraph{}We first introduce some notation. Let 
$\underline{\mathbf{d}}=(\mathbf{d}^{(1)}, \ldots, \mathbf{d}^{(r)})$
be a sequence of dimension vectors for $\Delta$ and let us set 
$\mathbf{d}=\sum_t \mathbf{d}^{(t)}$. We denote by
$\mathcal{F}_{\underline{\mathbf{d}}}$ the  variety of flags
$$0 =L_0\subset L_1 \subset \cdots \subset L_r=\bigoplus_i k^{\mathbf{d}_i}$$
of $I$-graded vector spaces such that 
$dim\;L_t/L_{t-1}=\mathbf{d}^{(t)}$. This is a smooth projective 
variety.
Now consider
$$X_{\underline{\mathbf{d}}}=\{((L_t),x) \in 
\mathcal{F}_{\underline{\mathbf{d}}} \times E(\mathbf{d})\;|\;
x(L_i) \subset L_{i-1}  \}.$$
This is a vector bundle over $\mathcal{F}_{\underline{\mathbf{d}}}$ 
of rank $\sum_i\sum_{s<t}\mathbf{d}^{(s)}_i
\mathbf{d}^{(t)}_{i-1}$. In
particular, it is a smooth (quasi-projective) variety. On the other 
hand, the second projection
$\pi:X_{\underline{\mathbf{d}}} \to E(\mathbf{d})$ is a projective morphism.
\begin{prop} Let $\mathbf{m}$ be a multipartition and let $x \in 
\mathcal{O}_{\mathbf{m}}$. Set
\begin{equation}\label{E:1}
\mathbf{d}^{(t)}=dim\; \mathrm{Ker}\; (x^t)/\mathrm{Ker}\;(x^{t-1}).
\end{equation}
Then $\pi(X_{\underline{\mathbf{d}}})=
\overline{\mathcal{O}_{\mathbf{m}}}$ and $\pi: X_{\underline{\mathbf{d}}} \to 
\overline{\mathcal{O}_{\mathbf{m}}}$
is a resolution of singularities.\end{prop}
\noindent
\textit{Proof.} Since $\pi$ is projective and $G_{\d}$-equivariant, its image is a
closed $G_\d$-stable subvariety. Moreover,
since $X_{\underline{\d}}$ is irreducible, so is $\pi(X_{\underline{\d}})$. Therefore
$\pi(X_{\underline{\d}})=\overline{\mathcal{O}}$ for some orbit $\mathcal{O}$.
Now let $x \in \pi(X_{\underline{\d}})$. By definition, the points of the fiber of 
$\pi$ over
$x$ are the flags $(L_t) \in \mathcal{F}_{\underline{\d}}$ which are compatible with $x$. 
It is
easy to see that such a flag is unique if and only if $x \in \mathcal{O}_{\m}$. 
On the other hand,
if $(L_t)$ is a flag compatible with $x$ and if $(L_t) \in \pi^{-1}(y)$ with $y \in
\mathcal{O}_{\m}$ then $x + uy \in \mathcal{O}_{\m}$ for $u$ in a dense 
subset of $\C^*$. Thus $x \in \overline{\mathcal{O}_{\m}}$. The Proposition follows.
\qed

\vspace{.2in}

\paragraph{}Observe that when $\Delta=\Delta_1^{(1)}$ and 
$\mathcal{O}_{\mathbf{m}}$ is the regular nilpotent orbit
then $X_{\underline{\mathbf{d}}}=T^*\mathfrak{B}$ where $\mathfrak{B}$ is the 
flag variety of $GL_d$ and the above resolution
is just the Springer resolution of the nilpotent cone.
This construction of desingularizations of orbit 
closures has been generalized to arbitrary finite
type quivers by Reineke (\cite{Re1}).

\begin{cor} The dimension of the orbit $\mathcal{O}_{\m}$ is
$$dim\;\mathcal{O}_{\m}=\sum_i \sum_{s<t} \d_i^{(s)}(\d_i^{(t)}+\d_{i-1}^{(t)}).$$
\qed
\end{cor}

\section{Hall algebras and canonical bases}

In this section we let $\Delta$ stand for either $\Delta_n$ or $\Dtn$.

\vspace{.2in}

\paragraph{\textbf{2.1.}} We briefly recall the definition of 
Ringel-Hall algebras. Let $k=\mathbb{F}_q$ be a finite field
with $q$ elements and let $\mathcal{C}=Mod_k(\Delta)$.
It is known that $\mathcal{C}$ satisfies the following conditions
\begin{enumerate}
\item[i)] $dim\;\mathrm{Ext}^i(M,N) < \infty$ for any $M,N \in \mathcal{C}$, $i=0,1$,
\item[ii)] $\mathrm{Ext}^i(M,N)=0$ for any $i>1$ and $M,N \in \mathcal{C}$.
\end{enumerate}
Introduce two bilinear forms on $K_0(\mathcal{C}) \simeq \Z^I$ by the formulas
$$ \langle M,N \rangle = 
dim\;\text{Hom}(M,N)-dim\;\text{Ext}^1(M,N),$$
$$ (M,N)=\langle M,N 
\rangle + \langle N,M\rangle.$$
Let us consider a $\C$-vector space $\U_{\mathcal{C}}$ (or simply $\U$) with a basis 
$\{[M]\}$ indexed by the isomorphism
classes of objects in $\mathcal{C}$. Let us choose $v \in \C$ such 
that $v^{-2}=q$. Finally let us denote by $\mathbf{d}(M) \in \N^I$ the dimension
of a representation $M$.

\paragraph{}We define an algebra structure on $\U$ by setting
\begin{equation}\label{E:2}
[M] \cdot [N]=v^{\d(M)\cdot \d(N)-\langle N,M\rangle}\sum_{[P]} F_{M,N}^P [P]
\end{equation}
where 
$$\d_1\cdot\d_2=dim\;(G_{\d_1+\d_2}/G_{\d_1} \times G_{\d_2})=2\sum_i\d_{1,i}
\d_{2,i},$$
$$F_{M,N}^P=\#\{Q \subset P\;|\; Q \simeq M, P/Q  \simeq N\}.$$
In order to define a bialgebra structure, we need to slightly enlarge 
$\U$. Let $\widetilde{\U}$ be the algebra
generated by $\U$ and elements $K_M$, $M \in K_0(\mathcal{C})$ subject to
the relations $K_M K_N=K_{M+N}$ and $K_M [P]=v^{-(M,P)} [P]K_M$. 
We define a bialgebra structure on
$\widetilde{\U}$ by putting
$$\Delta(K_M)=K_M \otimes K_M,$$
\begin{equation}
 \Delta([P])=\sum_{[M],[N]}
v^{-\d(M)\cdot\d(N)-\langle N,M\rangle}
\frac{|\text{Aut}(M)|\cdot|\text{Aut}(N)|}{|\text{Aut}(P)|} 
F^P_{M,N} [M]K_N \otimes [N].
\end{equation}
Finally, define a scalar product on $\widetilde{\U}$ by the formula
$$\langle 
[M]K_I,[N]K_J\rangle_G=v^{-(I,J)-\d(M)\cdot\d(N)}
\delta_{M,N}
\frac{(1-v^2)^{|\d(M)|}}{|\text{Aut}(M)|}$$
where we set $|\d|=\sum_i \d_i$.
\begin{theo}[\cite{Green}] $(\widetilde{\U}, \cdot, \Delta)$ is a bialgebra 
and $\langle\;,\; \rangle_G$ is a nondegenerate
Hopf pairing, i.e
$$\langle xy,z \rangle_G=\langle x \otimes y, \Delta(z)\rangle_G.$$
\end{theo}

\vspace{.1in}

\paragraph{\textbf{Remarks.}} i) The algebra $\U$ is clearly 
$K_0(\mathcal{C})=\Z^I$-graded. The same holds for $\widetilde{\U}$
where the elements $K_N$ are put in degree $0$.\\
ii) Our conventions differ slightly from the usual
ones regarding Ringel-Hall algebras. The standard basis $(\{M\})$ of $\U$ is related to 
ours by the change of variables $[M]=v^{|\d(M)|-\frac{1}{2}\d(M)\cdot \d(M)}\{M\}$,
and we renormalized the scalar
product $\langle\,,\,\rangle_G$ by the factor $(1-v^2)^{|\d(M)|}$. This will
make the relation to symmetric functions more transparent.

\vspace{.2in}

\paragraph{\textbf{2.2.}} Recall that the set of isomorphism classes of objects
in $Mod_k(\Delta)$ is in bijection with a set of multipartitions. 

\begin{prop}[\cite{Mac}, \cite{R3}] Let $\m,\mathbf{n},\mathbf{l}$ be arbitrary multipartitions. 
There
exists a polynomial $P_{\m,\mathbf{n}}^{\mathbf{l}} \in \Z[t]$
  such that for any finite field $k$ with $q$ elements we have
$F_{\m,\mathbf{n}}^{\mathbf{l}}=P_{\m,\mathbf{n}}^{\mathbf{l}}(q)$.
\end{prop}

The nice consequence of the above Proposition is that we now may
(and we now will) consider the Hall algebras 
$\U_{\mathcal{C}}$ and
$\widetilde{\U}_{\mathcal{C}}$ for a \textit{generic} $q$, 
i.e as $\C[v,v^{-1}]$-algebras. When $\Delta=\Delta_1^{(1)}$ then
the isomorphism classes of objects are labelled by usual partitions and the above
polynomials are the \textit{Hall polynomials} (see \cite{Mac}).

\vspace{.2in}

\paragraph{\textbf{2.3.}} Let us denote by $\C_{G_{\mathbf{d}}}
[E(\mathbf{d})]$ the set of $G_{\mathbf{d}}$-invariant functions on $E(\mathbf{d})$.
Observe that we have a natural identification of vector spaces
$$\U \simeq \bigoplus_{\mathbf{d}} \C_{G_{\mathbf{d}}}[E(\mathbf{d})], \qquad
[\mathbf{m}] \mapsto \mathbf{1}_{\mathcal{O}_{\mathbf{m}}}.$$
In this context, Lusztig gave in \cite{Lusbook} a geometric construction of the generic
Hall algebra $\U$, in which
one works over the algebraically closed field $\overline{\mathbb{F}_q}$ and where
$G_{\mathbf{d}}$-invariant functions are replaced by (complexes of) 
$G_{\mathbf{d}}$-equivariant constructible sheaves (see Saito's paper in this 
volume). Observe that in our case there are only finitely many $G_{\d}$-orbits
and these induce an algebraic stratification of $E(\d)$.
In particular, for every orbit $\mathcal{O}_{\m} \subset E(\d)$, let 
$IC_{\mathcal{O}_{\mathbf{m}}}$ be the simple
$G_{\mathbf{d}}$-equivariant perverse sheaf associated to the constant sheaf
${\C}_{\mathcal{O}_{\m}}$. We set
$$\mathbf{f}_{\m}=v^{dim\;\mathcal{O}_{\m}} \mathbf{1}_{\mathcal{O}_{\mathbf{m}}},$$
$$\mathbf{b}_{\m}=\sum_{i,\m'} v^{-i + dim\;\mathcal{O}_{\m}-dim\;
\mathcal{O}_{\m'}}
 dim\;\mathcal{H}^i_{\mathcal{O}_{\m'}}(IC_{\mathcal{O}_{\m}})\f_{\m'},$$
where $\mathcal{H}^i_{\mathcal{O}_{\m'}}(IC_{\mathcal{O}_{\m}})$ is the stalk
over any point of $\mathcal{O}_{\m'}$ of the $i$th cohomology sheaf of 
$IC_{\mathcal{O}_{\m}})$. The bases $\{\f_{\m}\}_{\m}$ and $\{\b_{\m}\}_{\m}$ are
called the \textit{PBW basis} and the \textit{canonical basis} respectively.
Observe that Green's scalar product restricts to a nondegenerate form on $\U$.
It will be important for us to consider the \textit{dual PBW basis} $\{\f^*_{\m}\}$
and the \textit{dual canonical basis} $\{\b^*_{\m}\}$ with respect to $\langle\;,\;
\rangle_G $. Note that $\f^*_{\m}$ and $\f_{\m}$ only differ by a scalar in 
$\C[v,v^{-1}]$.
For $\d \in \N^I$ we simply denote by $\f_{\d} (=\b_{\d})$ the basis element 
associated to the unique closed (=semisimple) orbit. 

\vspace{.2in}

\paragraph{\textbf{2.4.}} Finally, we note the following well-known recursive 
procedure to compute the canonical basis $\{\b_{\m}\}$ of $\U$. The Verdier
duality functor on the (derived) category of constructible $G_{\d}$-sheaves
on $E(\d)$ descends to a semilinear involution $u \mapsto \overline{u}$ on
$\U$. This involution can be computed~:

\begin{prop} The involution $u \mapsto \overline{u}$ is the unique ring homomorphism
$\U \to \U$ such that $\overline{v}=v^{-1}$ and $\overline{\f_{\d}}=\f_{\d}$ for
all $\d \in \N^I$.
\end{prop}
\noindent
\textit{Proof.} It is standard that the involution $u \mapsto \overline{u}$ has the
above properties. Let $\U' \subset \U$ be the subalgebra generated 
$\f_{\d}$ for $\d \in \N^I$. It is enough to prove that $\U'=\U$. We argue by 
induction on the (partial) degeneration order on the set of orbits. Let
$\m$ be a multipartition and assume that $\f_{\m'} \in \U'$ for all $\m' \neq \m$ such
that $\mathcal{O}_{\m'} \subset \overline{\mathcal{O}_{\m}}$. Pick $x \in 
\mathcal{O}_{\m}$ and define dimension
vectors $\d^{(t)}$, $t=1, \ldots, r$ as in (\ref{E:1}). Then from Proposition~1.1. we
have
$$\f^{\d_{r}} \cdots \f^{\d_{1}} \in v^{\Z}\f_{\m} \oplus 
\bigoplus_{\m'<\m} \C[v,v^{-1}] \f_{\m'}$$
and hence $\f_{\m} \in \U'$. \qed

\vspace{.2in}

The properties of simple perverse sheaves now translate as follows.
\begin{prop} The canonical basis element $\b_{\m}$ is characterized by the two
conditions
\begin{enumerate}
\item[i)] $\overline{\b_{\m}}=\b_{\m}$,
\item[ii)] $\b_{\m} \in \f_{\m} \oplus \bigoplus_{\m' < \m} v\C[v] \f_{\m'}.$
\end{enumerate}
\end{prop}
\qed

\section{Hall algebras and symmetric functions}

\paragraph{\textbf{3.1.}}In this section we assume that 
$\Delta= \Delta^{(1)}_1$ and set $\U^{(1)}_1=\U_{Mod_k(\Delta)}$. In that case, the set
of isomorphism classes of objects in $Mod_k(\Delta)$ is in 
natural bijection with the set $\Pi$ of all
partitions, and we will denote by the same letter $\lambda$ the 
isoclass corresponding to a partition $\lambda$.

\vspace{.2in}

\paragraph{\textbf{3.2.}} Let $\Gamma=\C[x_1, x_2, 
\ldots ]^{\mathfrak{S}_{\infty}}$ be Macdonald's ring of symmetric polynomials
in infinitely many variables. This ring is equipped with various bases. Following
classical notations, we denote by $p_{\lambda}$ the power-sum symmetric function,
by $e_{\lambda}$ the elementary symmetric function, and by $s_{\lambda}$ the Schur
function. 

\begin{theo}[\cite{Hall}] The assignement $\f_{(1^r)} \mapsto e_r$ extends to an 
isomorphism of algebras $\phi:\U^{(1)}_1 \stackrel{\sim}{\to} 
\Gamma \otimes_{\C} \C[v,v^{-1}]$. Moreover, $\phi(\f_{\lambda})=P_{\lambda}(v^2)$
where $P_\lambda$ is the Hall-Littlewood polynomial.\end{theo}

Now let us consider the scalar product $\langle\;,\;\rangle_G$. From the above theorem
one can deduce by direct computation~:
\begin{prop} The scalar product $\langle \;,\;\rangle=\phi^*\langle\;,\;\rangle_G$ 
coincides with the Hall-Littlewood scalar product~:
$$\langle p_{\lambda},p_{\mu} \rangle=\delta_{\lambda,\mu} z_{\lambda}
\frac{1}{\prod_i (1-v^{2\lambda_i})},$$
where $z_{(1^{m_1}2^{m_2}\cdots)}=\prod_i m_i!i^{m_i}$.
\end{prop}

\vspace{.1in}

Observe that the form $\langle\;,\;\rangle$ on $K_0(\mathcal{C})$ vanishes here.
It follows that $\U_1^{(1)}$ can itself be equipped with a coproduct by the formula
$$\Delta_0([P])=\sum_{[M],[N]}v^{-\d(M)\cdot\d(N)}
\frac{|\text{Aut}(M)|\cdot|\text{Aut}(N)|}{|\text{Aut}(P)|} 
F^P_{M,N} [M] \otimes [N].$$
From Green's theorem (Theorem~2.1) and the above proposition we see that $\phi:
\U^{(1)}_1 \stackrel{\sim}{\to} \Gamma\otimes_\C \C[v,v^{-1}]$ is an isomorphism of
bialgebras, where $\Gamma \otimes_{\C} \C[v,v^{-1}]$ is equipped with Zelevinsky's
coproduct (see \cite{Zel2}).

\vspace{.1in}

Finally, the canonical basis also admits a nice interpretation in terms of symmetric 
functions~:
\begin{prop}[\cite{Lus1}] We have $\phi(\b_{\lambda})=s_{\lambda}$. Hence
$\phi(\b^*_{\lambda})=S_{\lambda}(v)$ where $S_{\lambda}(v)$ is the dual
Schur function. \end{prop}

\vspace{.1in}

We will give a geometric interpretation of the above theorem in Section~7.

\section{Hall algebras and quantum groups}

\paragraph{\textbf{4.1.}}Let $A=(a_{ij})_{i=1,\ldots, r}$ be
a symmetric generalized Cartan matrix. Associated to $A$ is a complex
Kac-Moody algebra $\g$ with a Cartan decomposition $\g=\n_- \oplus \h \oplus \n_+$.
The quantized enveloping algebra $\U_v(\g)$ of $\g$ is the $\C(v)$-algebra generated by
elements $E_i, F_i, K_i^{\pm 1}$ for $i=1, \ldots, r$ satisfying the relations
$$K_iK_i^{-1}=K_i^{-1}K_i=1, \qquad K_iK_j=K_jK_i,$$
$$K_iE_j=v^{a_{ij}}E_jK_i, \qquad K_iF_j=v^{-a_{ij}}F_jK_i, \qquad [E_i,F_j]=
\delta_{ij}\frac{K_i-K_i^{-1}}{v-v^{-1}},$$
$$\sum_{k=0}^{1-a_{ij}} (-1)^k E_i^{(k)}E_jE_i^{(1-a_{ij}-k)}=\sum_{k=0}^{1-a_{ij}} 
(-1)^k F_i^{(k)}F_jF_i^{(1-a_{ij}-k)}=0\qquad\text{if}\;i\neq j,$$
where as usual
$$[k]=\frac{v^k-v^{-k}}{v-v^{-1}}, \qquad [k]!=[1]\cdots [k], \qquad
E_i^{(k)}=\frac{E_i^k}{[k]!}, \qquad F_i^{(k)}=\frac{F_i^k}{[k]!}.$$
The Lusztig integral form $\U_{A}(\g)$ is the $\C[v,v^{-1}]$-subalgebra generated
by $E_i^{(k)}, F_j^{(k)}$ and $K_i^{\pm 1}$ for $i =1, \ldots, r$ and $k \geq 1$. 
The above
algebras are equipped with a Hopf algebra structure with a coproduct
$$\Delta(K_i)=K_i \otimes K_i, \qquad \Delta(E_i)=E_i \otimes K_i + 1 \otimes E_i,
\qquad \Delta(F_i)=F_i \otimes 1 + K_i^{-1} \otimes F_i.$$
We will denote by $\U_A^-(\g) \subset \U_A(\g)$ the subalgebra generated by
$F^{(k)}_i$ for $i=1, \ldots, r$ and $k \geq 1$.

\vspace{.1in}

\paragraph{} To an arbitrary quiver $Q$ with vertex set $I$ and edge set $\Omega$
is associated a Cartan matrix $A$ as above~: we put $r=|I|$ and let $a_{ij}$ be the
number of (unoriented) edges between $i$ and $j$. In particular, $\Delta_n$
corresponds to the Lie algebra $\mathfrak{sl}_{n+1}$ while $\Dtn$ corresponds
to the affine Lie algebra $\widehat{\mathfrak{sl}}_n$ (for $n >1$).

\vspace{.2in}

\paragraph{\textbf{4.2.}} Let us first assume that $\Delta=\Delta_n$ and set
$\U_n=\U_{Mod_k(\Delta_n)}$. In that case
we have the following theorem of Ringel~:
\begin{theo}[\cite{R}] The assignment $F^{(k)}_i \mapsto [S_i^{\oplus k}]$ for 
$i=1, \ldots, n$
and $k \geq 1$
extends to an isomorphism of algebras $\phi:\U_A^-(\mathfrak{sl}_{n+1}) 
\stackrel{\sim}{\to} \U_n$.
\end{theo}

\vspace{.2in}

\paragraph{\textbf{4.3.}} For the rest of this section we assume that
$\Delta=\Dtn$ and set $\U^{(1)}_n=\U_{Mod_k(\Dtn)}$. The map 
$\phi$ defined above gives rise to an embedding
$\U^-_A(\widehat{\mathfrak{sl}}_n) \to \U^{(1)}_n$ but it is no longer surjective
(this can easily be checked by comparing graded dimensions).
Let us denote by $\mathbf{C}_n \subset \U^{(1)}_n$ the image of $\phi$. We first
give a description of $\mathbf{C}_n$.

\vspace{.1in}

\paragraph{}Recall that the set of isomorphism classes of objects in $\mathcal{C}$ is
indexed by (ordered) $n$-tuples of partitions $\mathbf{m}=(\lambda^{(1)}, \ldots,
\lambda^{(n)})$. Call such a multipartition \textit{aperiodic} if, for every $k >0$
there exists a least one $i$ such that $k$ is not a part of $\lambda^{(i)}$. Let us
denote by $\mathcal{M}$ the set of all multipartitions and by $\mathcal{M}^{ap}$ the
subset of aperiodic multipartitions.

\begin{theo}[\cite{Lus3.5}] We have 
$$\mathbf{C}_n=\bigoplus_{\m \in \mathcal{M}^{ap}} \C[v,v^{-1}] \b_{\m}.$$
\end{theo}

\vspace{.2in}

\paragraph{\textbf{4.4.}} Let us call a multipartition $\m=(\lambda^{(1)}, \ldots,
\lambda^{(n)})$ \textit{completely periodic} if we have $\lambda^{(1)}=\ldots
=\lambda^{(n)}$, and let $\mathcal{M}^{per} \subset \mathcal{M}$ be the set
of all such multipartitions. Consider the space
 $$\mathbf{K}=\bigoplus_{\m \in \mathcal{M}^{per}} \C[v,v^{-1}] \b^*_{\m}.$$
Observe that $\mathbf{K}$ has nonzero terms only in periodic degrees
$\mathbf{d}=(d, \ldots,d)$ (and hence it is in fact $\Z$-graded).

\begin{theo}[\cite{S1}] The following holds~:
\begin{enumerate}
\item[i)] $\mathbf{K}$ is a central subalgebra of $\U^{(1)}_n$,
\item[ii)] The multiplication map induces an isomorphism of algebras
$$\mathbf{K} \otimes_{\C[v,v^{-1}]} \mathbf{C}_n \to \U^{(1)}_n,$$
\item[iii)] The algebra $\mathbf{K}$ is stable under the involution $u \mapsto \overline{u}$ 
and, as a $\Z$-graded algebra,
$\mathbf{K} \simeq \C[v,v^{-1}][z_1,z_2, \ldots]$ where $deg\;(z_l)=l$ for $l \geq 1$.
\end{enumerate}
\end{theo}

Observe that if $n=1$ then $\mathbf{C}_n$ is not defined and $\U^{(1)}_1=\mathbf{K}$. Part iii)
of the above theorem can then be deduced from Theorem~3.1 and the well-known fact that 
$\Gamma \simeq 
\C[z_1, z_2, \ldots]$ as a $\Z$-graded algebra. However, we will see that, as in the
case $n=1$, there is \textit{for every $n$} a geometrically defined \textit{canonical} 
identification $\mathbf{K} \simeq \Gamma \otimes \C[v,v^{-1}]$.

\vspace{.1in}

\paragraph{\textbf{Remarks.}} i) Theorem~4.3 can be reinterpreted as follows~: 
$\U^{(1)}_n \simeq
\U_A^-(\widehat{\mathfrak{gl}}_n)$ where $\widehat{\mathfrak{gl}}_n$ is defined as
the (universal) central extension
$$0 \to \C c \to \widehat{\mathfrak{gl}}_n \to \mathfrak{gl}_n[z,z^{-1}] \to 0.$$
The algebra $\mathbf{K}$ is then the quantum analogue of the algebra generated by
$xz^k$ for $k \geq 1$, where $x \in Z(\mathfrak{gl}_n)$.

ii) The algebra $\mathbf{K}$ admits another description. Let
$\mathbf{e}'_i
\in \mathrm{End}(\U^{(1)}_n)$ be the adjoint (under $\langle \;,\;\rangle_G$) of the left
multiplication by $F_i$ (the Kashiwara operator familiar in crystal graph theory, see e.g.
\cite{HK}). 
Then $\mathbf{K}=\cap_i \mathrm{Ker}\;\mathbf{e}'_i$. In particular, $\mathbf{K}$
corresponds to the ``highest weight'' subspace (of the crystal graph of $\U^{(1)}_n$)
with respect to the action of the quantum group $\U_A(\widehat{\mathfrak{sl}}_n)$.
The above results are in accordance with (and were partly motivated by) \cite{LTV} where
the crystal graph of $\U^{(1)}_n$ was completely described, and found to be a disjoint union
of copies of the crystal graph of $\U^-_A(\widehat{\mathfrak{sl}}_n)$ indexed
by partitions.

iii) In the limit $n \to \infty$ the diagrams $A^{(1)}_n$ and $A_n$ tend to the common 
limit
$A^{\infty}_{\infty}$ (i.e. the quiver with vertex set $I+\Z$ and edge set $\Omega=\{
(i,i+1)\}$). In that situation, $\mathbf{U}^{(1)}_{\infty}=\mathbf{U}_{\infty}=
\mathbf{C}_{\infty}=\U^-_A(\mathfrak{sl}_{\infty})$.

iv) Ringel also studied a ``degenerate'' version of the Hall
algebra where one considers representations of the given quiver over $\C$ and uses an
Euler characteristic analogue of the convolution product (\ref{E:2}). This yields a
construction of the \textit{Lie algebra} of $\mathfrak{sl}_n$ or
$\widehat{\mathfrak{gl}}_n$ rather than of the quantum group (see \cite{R2}).

\section{Flag varieties and Hecke algebras}

In this section we (brielfly) introduce and fix some notation for
the more classical geometric objects
in type $A$ representation theory such as (partial, affine) flag varieties, (affine) 
Hecke algebras and some related algebras. Good references for these are \cite{BG}
or \cite{CG}. 

\vspace{.2in}

\paragraph{\textbf{5.1.}} Let $\mathfrak{S}_r$ denote the symmetric group on $r$ letters.
We denote by $s_1, \ldots, s_{r-1}$ the simple reflections, so that $\mathfrak{S}_r$ may be
presented as the group generated by $s_i, i=1, \ldots, r-1$ with relations
\begin{equation}\label{E:4}
s_i^2=1, \qquad s_is_j=s_js_i \text{\;for\;}i \neq j \pm 1,
\end{equation}
\begin{equation}\label{E:5}
s_is_{i+1}s_i=s_{i+1}s_is_{i+1}.
\end{equation}
The extended affine symmetric group is $\widehat{\mathfrak{S}}_r=\mathfrak{S}_r
\ltimes \Z^r$. It is isomorphic to the group generated by $s_i$ for $i \in \Z/r\Z$ and $\tau$
subject to the relations (\ref{E:4}), (\ref{E:5}), $\tau^r=1$ and
$\tau s_i =s_{i+1}\tau$.\\
As usual we will denote by $l: \mathfrak{S}_r \to \N$ and
$l: \widehat{\mathfrak{S}}_r \to \N$ the length functions. By definition, we have
$l(s_i)=1$ and 
$l(\tau)=0$.

\vspace{.1in}

\paragraph{} Now fix $\mathbf{d}=(d_1, \ldots, d_n) \in \N^n$ such that $\sum_i d_i=r$.
Associated to $\mathbf{d}$ are parabolic subgroups $\mathfrak{S}_{\mathbf{d}} \subset
\mathfrak{S}_r$ and $\widehat{\mathfrak{S}}_{\mathbf{d}} \subset \widehat{\mathfrak{S}}_r$
generated by $s_i$ for $i \not\in \{d_1, d_1+d_2, \cdots, d_1 + \cdots + d_n\}$. Observe 
that $\widehat{\mathfrak{S}}_{\d}$ is always finite. 

\vspace{.2in}

\paragraph{\textbf{5.2.}} Let $k$ be either $\mathbb{F}_q$ or $\overline{\mathbb{F}_q}$
Let $L$ be an $r$-dimensional $k$-vector space. We denote
by $\mathcal{B}$ the variety of all complete flags 
$$0=L_0 \subset L_1 \cdots \subset L_r=L, \qquad dim\;L_i=i$$
(a smooth projective variety). The group $GL(L)$ acts transitively on $\mathcal{B}$.
The set of $GL(L)$-orbits on $\mathcal{B} \times \mathcal{B}$ is in bijection with the
symmetric group $\mathfrak{S}_r$. For $w \in \mathfrak{S}_r$ we denote by $\mathcal{O}_w$
the corresponding orbit. 

\vspace{.1in}

\paragraph{}Alternatively, we may fix one flag $(L_i)_i \in \mathcal{B}$ and 
denote by $B \subset GL(L)$ its stabilizer (a Borel subgroup). This provides us with an 
identification
$\mathcal{B} \simeq GL(L)/B$. The set of $B$-orbits on $\mathcal{B}$ is also in natural
bijection with $\mathfrak{S}_r$ and we denote by $\mathcal{B}_w$ the orbit corresponding to
$w$ (a Schubert cell). It is known that $dim\;\mathcal{B}_w=l(w)$.

\vspace{.1in}

\paragraph{}Now let $n \in \N$ and fix $\mathbf{d}=(d_1, \ldots, d_n) \in \N^n$ such that
$\sum_i d_i=r$. Let $\mathcal{F}_{\d}$ be the variety of $n$-step
partial flags
$$0=L_0 \subset L_1 \subset \cdots \subset L_n =L, \qquad dim\;L_i/L_{i-1}=d_i.$$
This is again a smooth projective variety.
As above, let $(L_i)_i$ be a fixed flag in $\mathcal{F}_{\d}$ and let $P \subset
GL(L)$ be its stabilizer (a parabolic subgroup). This identifies $\mathcal{F}_{\d}$ with
$GL(L)/P$. The $P$-orbits on $\mathcal{F}_{\d}$ are now parametrized by double
cosets $\mathfrak{S}_\d \backslash \mathfrak{S}_r/
\mathfrak{S}_\d$.
We again denote by $\mathcal{F}_{\d,w}$ the orbit corresponding to $w \in \mathfrak{S}_\d 
\backslash\mathfrak{S}_r/\mathfrak{S}_\d$ (a generalized Schubert cell).

\vspace{.1in}

\paragraph{}The Zariski closures $\overline{\mathcal{B}_w}$ or $\overline{\mathcal{F}_{\d,w}}$
of Schubert cells are called \textit{Schubert varieties}.

\vspace{.2in}

\paragraph{\textbf{5.3.}} We consider an affine analogue of the above notions. Set
$\mathbb{L}=k^r((z))$. By definition, a \textit{lattice} in $\mathbb{L}$ is a free
$k[[z]]$-submodule $L \subset \mathbb{L}$ of rank $r$. Let $\mathcal{G}$ denote the 
variety of all lattices. This infinite-dimensional
space is a countable union of (finite-dimensional) smooth projective varieties (
indeed, let us fix one lattice $L_0$; then for every $l \in \N$ let $X_l \subset 
\mathcal{G}$ be the smooth projective, finite-dimensional variety of lattices
$L$ with $z^lL_0 \subset L \subset z^{-l}L_0$; we have $X_l \subset X_{l+1}$ and
$\mathcal{G}=\cup_l X_l$.).

\vspace{.1in}

\paragraph{}Now let $\mathcal{X}$ be the variety of all infinite flags of 
lattices $ \cdots \subset L_i \subset L_{i+1} \subset \cdots $
such that $dim_{k} L_{i+1}/L_{i}=1$ and which
satisfy the periodicity condition $L_{i+r}=z^{-1}L_i$. This is again a union of 
finite-dimensional smooth projective varieties. The map $\mathcal{X} \to \mathcal{G},
(L_i) \mapsto L_0$ is a smooth locally trivial $B$-fibration. The group $GL_{k((z))}
(\mathbb{L})$ acts transitively on $\mathcal{X}$ and the set of 
$GL_{k((z))}(\mathbb{L})$-orbits on $\mathcal{X} \times \mathcal{X}$ is in natural
bijection with $\widehat{\mathfrak{S}}_r$ : for $j \in \Z/r\Z$, $s_j$ corresponds to
the orbit consisting of pairs $((L_i),(L'_i))$ such that $L_i=L'_i$ for all
$i \not\equiv j\;(\text{mod}\;r)$; the orbit $\{((L_i),(L'_i))\;|L_i=L'_{i+1} 
\;\forall\;i\}$ corresponds to $\tau$. For any $w \in \widehat{\mathfrak{S}}_r$ we
will write $\mathcal{O}_w$ for the associated orbit in $\mathcal{X} \times \mathcal{X}$.

\vspace{.1in}

As in Section~5.2, we may also fix a flag $(L_i) \in \mathcal{X}$ and denote by
$I \subset GL_{k((z))}(\mathbb{L})$ its stabilizer (an Iwahori subgroup). This identifies
$\mathcal{X}$ with $GL_{k((z))}(\mathbb{L})/I$. The set of $I$-orbits on $\mathcal{X}$
is also in natural bijection with $\widehat{\mathfrak{S}}_r$ and we will denote by
$\mathcal{X}_w$ the orbits corresponding to $w \in \widehat{\mathfrak{S}}_r$ (an affine
Schubert cell). As in the finite type case, we have $dim\;\mathcal{X}_w=l(w)$.

\vspace{.1in}

\paragraph{}Finally, let $\d=(d_1, \ldots, d_n) \in \N^n$ be such that $\sum_i d_i=r$. We
set
$$\mathcal{Y}_{\mathbf{d}}=\{(L_i)\;|L_i \subset L_{i+1}, \quad dim_k\;L_i/L_{i-1}=
d_{i\;\text{mod}\;n},\quad L_{i+n}=z^{-1}L_i\}.$$
As before, there is a (noncanonical) indentification $\mathcal{Y}_{\d}
\simeq GL_{k((z))}(\mathbb{L})/P$ for some (parahoric) subgroup $P$; the $P$-orbits
$\mathcal{Y}_w$ on $\mathcal{Y}_{\d}$ are indexed by elements $w \in 
\widehat{\mathfrak{S}}_{\d}
\backslash \widehat{\mathfrak{S}}_r / \widehat{\mathfrak{S}}_\d$ (see \cite{Lus3} 
for a formula
for $dim\; \mathcal{Y}_w$).

\vspace{.2in}

\paragraph{\textbf{5.4.}} Define the Hecke algebra $\mathbf{H}_r$ (resp. the affine 
Hecke algebra $\widehat{\mathbf{H}}_r$) as the unital $\C[v,v^{-1}]$-algebra generated
by elements $T_w$, $w \in \mathfrak{S}_r$ (resp. $w \in \widehat{\mathfrak{S}}_r$) 
satisfying the following relations
$$T_wT_{w'}=T_{ww'} \qquad \text{if}\; l(ww')=l(w)+l(w'),$$
$$(T_{s_i}+1)(T_{s_i}-v^{-2})=0,\qquad \text{for\;all}\;i.$$
We set $\tilde{T}_{w}=v^{l(w)}T_w$ for every $w\in 
\mathfrak{S}_r$ (resp. $w \in \widehat{\mathfrak{S}}_r$).

In the case of an affine Hecke algebra, there is another very useful presentation, due
to Bernstein.

\begin{theo}[see \cite{Lus4}] The algebra $\widehat{\mathbf{H}}_r$ is isomorphic to the 
algebra generated by $T^{\pm 1}_i$, $i=1, \ldots, r-1$ and $X_j^{\pm 1}$, $j =1, \ldots, r$
with relations
\begin{alignat*}{2}
&T_i\,T_i^{-1}=1=T_i^{-1}\,T_i,\qquad & \qquad &(T_i+1)(T_i-v^{-2})=0,\\
&T_i\,T_{i+1}\,T_i=T_{i+1}\,T_i\,T_{i+1},\qquad&\qquad&
|i-j|>1\Rightarrow T_i\,T_j=T_j\,T_i,\\
&X_i\,X_i^{-1}=1=X_i^{-1}\,X_i,\qquad&\qquad &X_i\,X_j=X_j\,X_i,\qquad\\
&T_i\,X_i\,T_i=v^{-2}X_{i+1},\qquad&\qquad 
&j\not= i,i+1\Rightarrow X_j\,T_i=T_i\,X_j.
\end{alignat*}
\end{theo}

The isomorphism between the two presentations is such that
$T_{s_i}\mapsto T_i$ and $\tilde{T}_\lambda^{-1}\mapsto X_1^{\lambda_1}\cdots
X_r^{\lambda_r}$ if $\lambda=(\lambda_1,\ldots,\lambda_r)\in \Z^r$ is
\textit{dominant}, i.e if $\lambda_1 \geq \lambda_2 \geq \cdots \geq \lambda_r$. 
The center of $\widehat{\mathbf{H}}_r$ is $Z(\widehat{\mathbf{H}}_r)=
\C[v,v^{-1}][X_1^{\pm 1},\ldots ,X_r^{\pm 1}]^{\mathfrak{S}_r}$. 

\vspace{.2in}

\paragraph{}The relevance of Hecke algebras for us and the link with the previous paragraphs
comes from the following geometric
construction. Let us assume that $k=\mathbb{F}_q$. Let us
set $G=GL_{k((z))}(\mathbb{L})$ and denote by $\C_G[\mathcal{X} \times
\mathcal{X}]$ the space of $G$-invariant functions whose supports lie on a finite union
of orbits. We equip $\C_G[\mathcal{X} \times \mathcal{X}]$ with the algebra structure
given by the convolution product
$$ f \star g = p_{13!}\left(p_{12}^*(f) \cdot p_{23}^*(g)\right)$$
where $p_{ij}: \mathcal{X} \times \mathcal{X} \times \mathcal{X} \to \mathcal{X} \times 
\mathcal{X} $ is the projection on the $i$th and $j$th components. Observe that this
product is well-defined since for any two orbits $\mathcal{O}_w$, $\mathcal{O}_{w'}$,
the projection $p_{13}: (\mathcal{O}_w \times \mathcal{X}) \cap (\mathcal{X} \times 
\mathcal{O}_{w'}) \to \mathcal{X} \times \mathcal{X}$ is proper.

\vspace{.2in}
In the next theorem, we specialize the Hecke algebra by setting $v^{-2}=q$.

\begin{theo}[\cite{IM}] The assignement $T_w \mapsto \mathbf{1}_{\mathcal{O}_w}$
defines an isomorphism of algebras $\widehat{\mathbf{H}}_{r} \simeq \C_G[\mathcal{X}
\times \mathcal{X}]$.\end{theo}

\vspace{.2in}

\paragraph{} It follows from the above Theorem that the structure constants for the convolution
algebra $\C_G[\mathcal{X} \times \mathcal{X}]$ are polynomials in $q$. Hence we may (and we will)
consider the convolution algebra as a $\C[v,v^{-1}]$-algebra.

\vspace{.1in}

Naturally, the same result holds in the case of the usual flag variety $\mathcal{B}$, and
we obtain in this way a geometric construction of the finite Hecke algebra 
$\mathbf{H}_r$.
 
\vspace{.2in}

\paragraph{\textbf{5.5.}} Let us now consider parabolic analogues of the constructions
in the previous paragraph. We again assume that $k=\mathbb{F}_q$.
Let us set
$$\mathcal{Y}=\sqcup_{\mathbf{d}} \mathcal{Y}_{\d},$$
where in the sum $\d$ runs over all sequences $(d_1, \ldots ,d_n) \in \N^n$
satisfying $\sum d_i=r$.

\vspace{.1in}

 The group $G=GL_{k((z))}(\mathbb{L})$ acts on $\mathcal{Y} \times 
\mathcal{Y}$. It is easy to check that the $G$-orbits in 
$\mathcal{Y}_\d \times \mathcal{Y}_{\d'}$ are indexed
by $\widehat{\mathfrak{S}}_\d \backslash \widehat{\mathfrak{S}}_r
/ \widehat{\mathfrak{S}}_{\d'}$. Without risk of confusion, we denote by $\mathcal{O}_w 
\subset \mathcal{Y}_\d \times \mathcal{Y}_{\d'}$ the orbit corresponding to $w$.

\vspace{.2in}

\paragraph{} For ${w} \in  \widehat{\mathfrak{S}}_\d \backslash 
\widehat{\mathfrak{S}}_r
/ \widehat{\mathfrak{S}}_{\d'}$, we set
$T_{{w}}=\sum_{\delta \in {w}} T_{\delta}$
and we let $\AH_{\mathbf{dd}'} \subset \AH_r$ be the $\C[v,v^{-1}]$-linear span of the
elements $T_{{w}}$ for ${w} \in \Si_{\mathbf{d}}
\backslash \Si_r/\Si_{\mathbf{d}'}$. Set $e_{\mathbf{d}}=\sum_{\delta
\in \Si_{\mathbf{d}}} T_{\delta}$. Then $\AH_{\mathbf{dd}'}=e_{\mathbf{d}} \AH_r
e_{\mathbf{d}'}$. Put
$$\widehat{\mathbf{S}}_{n,r}=\bigoplus_{\mathbf{d},\mathbf{d}'} \AH_{\mathbf{dd}'}.$$
We equip this space with the multiplication
$$e_{\mathbf{d}} h e_{\mathbf{d}'} \bullet e_{\mathbf{c}} h' e_{\mathbf{c}'}=
\delta_{\d'\mathbf{c}} e_{\mathbf{d}} h e_{\mathbf{d}'} h' e_{\mathbf{c}'} \in
\AH_{\mathbf{\d\mathbf{c}'}}\qquad {for\;all}\; h,h' \in \AH_r.$$
 
\vspace{.1in}

\paragraph{}Of course, replacing everwhere $\mathcal{Y}$ by $\mathcal{F}$ and 
$\widehat{\mathfrak{S}}_r$ by $\mathfrak{S}_r$, we obtain the definition of the algebra 
$\mathbf{S}_{n,r}$.

\vspace{.2in}

In the next theorem, we specialize the Hecke algebra by setting $v^{-2}=q$.

\begin{theo}[\cite{VV}]  The linear map $\widehat{\mathbf{S}}_{n,r} \to 
\C_{G}
[\mathcal{Y} \times \mathcal{Y}]$ such that $T_{{w}} \mapsto
\mathbf{1}_{\mathcal{O}_{{w}}}$
is an algebra isomorphism.
\end{theo}

Again a similar result holds in the finite type case. Finally, as in Section~5.4, we deduce
from the above Theorem that $\C_G[\mathcal{Y} \times \mathcal{Y}]$ is naturally a 
$\C[v,v^{-1}]$-algebra.

\vspace{.2in}

\paragraph{\textbf{5.6.}} As for Hall algebras, there is a construction of the Hecke
algebras involving (complexes of) equivariant constructible sheaves,
rather than $G$-invariant functions. We will not need such a construction and refer to 
e.g.  \cite{BG} or \cite{Vas} for details. A consequence
of this construction is the existence of canonical bases for these algebras, which we now
define. 

\vspace{.1in}

Let $w \in \widehat{\mathfrak{S}}_r$. We put
$$\mathbf{c}_w=\sum_{i,x} v^{-i+l(w)-l(x)}{dim}\;
\mathcal{H}^i_{\mathcal{X}_{x}}(IC_{\mathcal{X}_{w}}) \tilde{T}_x$$
where $IC_{\mathcal{X}_w}$ denotes the intersection cohomology complex
associated to $\C_{\mathcal{X}_{w}}$ and where $\mathcal{H}^i$ stands for
local cohomology. The elements $\{\mathbf{c}_w\;|\;w\in \widehat{\mathfrak{S}}_r\}$
form the canonical (or \textit{Kazhdan-Lusztig}) $\C[v,v^{-1}]$-basis of 
$\widehat{\mathbf{H}}_r$.

\vspace{.1in}

\paragraph{}Similarly, fix $\mathbf{d},\mathbf{d}'$
and let
${w} \in
\Si_{\mathbf{d}}\backslash {\Si}_r/\Si_{\mathbf{d}'}$. 
Let $L_{\mathbf{d}}$ be any point in $\mathcal{Y}_{\d}$ and
denote by $\mathcal{Y}_{{w},\mathbf{d}}$ the
fiber above $L_{\d}$ of the projection of $\mathcal{O}_{{w}} \to
\mathcal{Y}_{\d}$ on the first component. This is an algebraic
variety of dimension $y({w})$, say (an explicit formula for
$y({w})$ can be found in \cite{Lus3}). 
For every ${w} \in \Si_{\mathbf{d}}\backslash
{\Si}_r/\Si_{\mathbf{d}'}$ set
$$\mathbf{c}_{{w}}=\sum_{i,{x}}
v^{-i+y({w})-
y({x})} \;{dim}\;\mathcal{H}^i_{
\mathcal{Y}_{x,\mathbf{d}}}
(IC_{\mathcal{Y}_{w,\mathbf{d}}}) \tilde{T}_{x}.$$

The elements 
$\bigsqcup_{\mathbf{d},\d'}\{\mathbf{c}_w\;|\;w \in 
\Si_{\mathbf{d}}\backslash {\Si}_r/\Si_{\mathbf{d}'}\}$
form the canonical basis of $\widehat{\mathbf{S}}_{n,r}$.

\vspace{.2in}

\paragraph{}In addition, as for Hall algebras, there is a purely algebraic
characterization of the elements of the canonical bases of the algebras $\H_r$ and 
$\mathbf{S}_{n,r}$.

\vspace{.1in}

Let us first consider the semilinear involution $u \mapsto \overline{u}$ of 
$\widehat{\mathbf{H}}_r$ defined by $\overline{T_w}=T_{w^{-1}}^{-1}$.

\begin{theo}[\cite{KL}, \cite{KL2}]
For each $w \in {\Si}_r$ the element $\mathbf{c}_w
\in \AH_r$ is uniquely determined by
$$\mathrm{i)} \;\;\overline{\mathbf{c}_w}=\mathbf{c}_{w},\qquad \mathrm{ii)}\;\;
\mathbf{c}_w=\tilde{T}_w
+\sum_{x<w}P_{x,w}(v)\tilde{T}_x,\;\quad P_{x,w}(v)\in v\C[v].$$
\end{theo}

The polynomial $P_{x,w}(v)$ is the affine Kazhdan-Lusztig polynomial
of type $\tilde{A}_{r-1}$ associated to $w$ and $x$ 

\vspace{.1in}

Similarly, define a semilinear
involution $\tau : \AH_{\mathbf{dd}'} \to \AH_{\mathbf{dd}'}$ by $\tau (u)=
v^{-2l(w_{\mathbf{d}'})} \overline{u}$ where $w_{\mathbf{d}'} \in 
\widehat{\mathfrak{S}}_{\d'}$ is the longest element.
The previous theorem generalizes as follows.

\begin{theo} For each $w \in \Si_{\mathbf{d}}\backslash
{\Si}_r/\Si_{\mathbf{d}'}$
the element $\mathbf{c}_{w}$
is uniquely determined by the properties :
$$\mathrm{i)}\;\tau(\mathbf{c}_{w})=
\mathbf{c}_{w},\qquad
\mathrm{ii)}\;\mathbf{c}_{w}=\tilde{T}_{w}+
\sum_{x
<w} Q_{x,w}(v)
\tilde{T}_{x},\qquad
Q_{x,w}(v) \in v\C[v].$$
\end{theo}

\vspace{.1in}

\paragraph{}Again, we have given the results in the affine case but they
are all valid in the finite type case as well.

\section{Canonical bases and Kazhdan-Lusztig polynomials}

We now give the fundamental constructions (due to Lusztig and Zelevinsky) relating
representations of quivers of type $A$ to flag varieties. In particular, this will
identify the coefficients of the transition matrix between the canonical basis $\{\b_\m\}$
and the PBW basis $\{\f_{\m}\}$ with Kazhdan-Lusztig polynomials. The strategy is
as follows~: we construct a compactification of an orbit closure 
$\overline{\mathcal{O}_{\mathbf{m}}}$ which is a (generalized) Schubert variety, in such
a way that the stratification by Schubert cells restricts on 
$\overline{\mathcal{O}_{\mathbf{m}}}$ to the stratification
by $G_{\d}$-orbits and use Theorem~5.4. In this section, $k=\C$ or $k=\overline{\F}$.

\vspace{.2in}

\paragraph{\textbf{6.1.}} Let us first treat the finite type case. Let $\d=(d_1, \ldots
, d_n) \in \N^n$ be a dimension vector, and let us set $\sum d_i=r$ and
$$V_i=k^{d_i}, \qquad V^i=\bigoplus_{j \leq i} V_j, \qquad L=\bigoplus_i V_i.$$
Then $(V^i)$ is a partial flag belonging to $\mathcal{F}_{\d}$. Let us denote 
by $P\subset GL(L)$ its stabilizer. We construct a map $\iota:\; E(\d) \to 
\mathcal{F}_{\d}$ as follows. To a representation $(\phi_i) \in E(\d)$ we associate the 
flag defined by
$$U_i=\{(x_1, \ldots, x_n) \in V_1 \oplus \cdots \oplus V_n\;| x_{j+1}=\phi_j(x_j)\;
\text{for\;} j \geq i\}.$$

\begin{prop}[\cite{Zel}] The map $\iota$ is an embedding of $E(\d)$ as a dense open
subset of the (generalized) Schubert variety
$$\overline{\mathcal{F}_{w_0}}:=\{(U_i)\;|\; U_i \supset V^{i-1}\quad \text{for\;all\;}
i\}.$$
Moreover, each $G_{\d}$-orbit $\mathcal{O}_{\m} \subset E(\d)$ embeds
as a dense open subset of a $P$-orbit $\mathcal{F}_{w(\m)}$ for some $w(\m) \in
\mathfrak{S}_{\d} \backslash \mathfrak{S}_{r} / \mathfrak{S}_{\d}$. 
\end{prop}

\noindent
\textit{Proof.} It is clear that $\iota(E(\d)) \subset \overline{\mathcal{F}_{w_0}}$.
Conversely, it is easy to see that $(U_i) \in \overline{\mathcal{F}_{w_0}}$ belongs to
$\iota(E(\d))$ if and only if $U_i \cap \bigoplus_{j > i}V_j =0$ for all $i$ 
(an open condition).
Finally note that $G_{\d}$ embeds as a Levi subgroup of $P$ and that $\iota$ is 
$G_{\d}$-equivariant. \qed

\vspace{.2in}

\paragraph{\textbf{Remark.}} Let $\mathcal{U} \subset \overline{\mathcal{F}_{w_0}}$ be any
$P$-orbit and pick $(U_i) \in \mathcal{U}$. Then
$$(U_i/V^{i-1}, \; \phi_i=Id : U_i/V^{i-1} \to U_{i+1}/V^i)$$
defines a representation of $\Delta_n$. Let $\m$ be the multipartition describing this
representation. It is easy to check that $\iota(\mathcal{O}_{\m}) \subset \mathcal{U}$.
In particular, the assignement $\m \mapsto w(\m)$ is a bijection between the set of
$G_{\d}$-orbits in $E(\d)$ and the set of $P$-orbits (or Schubert cells)
in $\overline{\mathcal{F}_{w_0}}$.

\vspace{.2in}

\paragraph{}Now let us fix $\mathcal{O}_{\m} \subset E(\d)$. By the previous Proposition,
 the simple perverse sheaf $IC_{\mathcal{O}_{\m}}$ is the restriction to
$\overline{\mathcal{O}_{\m}}$ of the simple perverse sheaf 
$IC_{\overline{\mathcal{F}_{w_0}}}$. In particular, for any $\mathcal{O}_{\mathbf{n}}
\subset \overline{\mathcal{O}_{\m}}$ we have
$$dim\; \mathcal{H}^i_{\mathcal{O}_{\mathbf{n}}}(IC_{\mathcal{O}_{\m}})=
dim\; \mathcal{H}^i_{\mathcal{F}_{w(\mathbf{n})}}(IC_{\mathcal{F}_{w(\m)}}).$$
Next, consider the natural map $\pi: \mathcal{B} \simeq G/B \to G/P \simeq 
\mathcal{F}_{\d}$ given by
$$(0=U_0 \subset U_1 \cdots \subset U_{r}=L) \mapsto (U_{d_1} \subset U_{d_1+d_2} 
\subset \cdots U_{r}=L).$$
This is a locally trivial fibration with smooth projective fiber $P/B$. Thus, for any
$w_1,w_2 \in \mathfrak{S}_{\d} \backslash \mathfrak{S}_{r} / \mathfrak{S}_{\d}$ 
there holds
$$dim\; \mathcal{H}^i_{\mathcal{F}_{w_1}}(IC_{\mathcal{F}_{w_2}})=
dim\; \mathcal{H}^i_{\pi^{-1}(\mathcal{F}_{w_1})}(IC_{\pi^{-1}(\mathcal{F}_{w_2})}).$$
Finally, $\pi^{-1}(\mathcal{F}_w) =\sqcup_{y \in w} \mathcal{B}_y$ (a finite union of
$B$-orbits). Since $\mathcal{F}_w$ is irreducible, there exists a unique open $B$-orbit
$\mathcal{B}_{y_w}$ in $\pi^{-1}(\mathcal{F}_w)$ (corresponding to the longest
element $y(w) \in w$). Thus we obtain

\begin{equation*}
\begin{split}
dim\; \mathcal{H}^i_{\mathcal{O}_{\mathbf{n}}}(IC_{\mathcal{O}_{\m}})=&
dim\; \mathcal{H}^i_{\pi^{-1}(\mathcal{F}_{w(\mathbf{n})})}
(IC_{\pi^{-1}(\mathcal{F}_{w(\m)})})\\
=&dim\; \mathcal{H}^i_{\mathcal{B}_{y(w(\mathbf{n}))}}
(IC_{\mathcal{B}_{y(w(\m))}}).
\end{split}
\end{equation*}

\begin{cor} The coefficients of the transition matrix between the canonical
basis $\{\b_{\m}\}$ and the PBW basis $\{\f_{\m}\}$ of $\U_n$
are Kazhdan-Lusztig polynomials 
(of type $A_{\sum d_i -1})$, where $\d$ is determined by $\m$ via $\O(\m)\subset E_{\d}$.
\end{cor}

\vspace{.2in}

\paragraph{\textbf{6.2.}} Let us now briefly describe the affine analogue of the above
construction, which is due to Lusztig. Let $\d=(d_1, \ldots, d_n) \in \N^n$ and set 
$r=\sum d_i$. For clarity, we denote by $\overline{i} \in \Z/n\Z$ the class of an
integer $i \in \Z$. Set $V_i=k^{d_i}$, $L=\bigoplus_i V_i$ and $\mathbb{L}=L((z))$.
Let $ \cdots \subset\mathcal{L}_0 \subset \mathcal{L}_1 \subset \cdots $ be the 
$n$-periodic flag of lattices defined by
$$\mathcal{L}_0=L[[z]],$$
$$\mathcal{L}_i=\mathcal{L}_0 \oplus z^{-1} V_1 \oplus \cdots
 \oplus z^{-1}V_i\qquad \text{for\;}1 \leq i \leq n,$$
$$\mathcal{L}_{j+n}=z^{-1}\mathcal{L}_j\quad \text{for\;all}\;j.$$
 Let $P \subset GL_{k((z))}(\mathbb{L})$ be its stabilizer. We construct a map
$E(\d) \to \mathcal{Y}_{\d}$ as follows~: let $0 \leq \overline{i} <n$ be the residue of 
$i$ modulo $n$ and set $l(j)=(j-\overline{j})/n$;
let $(\phi_i) \in E(\d)$ and set
$$\tilde{\phi}_i=\begin{cases} \phi_i & \text{if}\; \overline{i} \neq 0\\
z^{-1}\phi_i& \text{if}\; \overline{i}=0,\end{cases}$$
and $\phi^h(x_i)=\tilde{\phi}_{i+h} \cdots \tilde{\phi}_{i+1} \tilde{\phi}_i(x_i)$
if $x_i \in V_i$. Note that for any fixed $x_i$, $\phi^h(x_i)=0$ for $h$ big enough.
Now we put $\iota((\phi_i))=(\mathcal{N}_j)$ where
$\mathcal{N}_j=\mathcal{L}_{j-1} \oplus z^{-l(j-1)-1} \mathcal{U}_{j-1}$
and
$$\mathcal{U}_j=\big\{\sum_{h \geq 0} \phi^h(x_{\overline{j}})\;| x_{\overline{j}}
\in V_{\overline{j}}\big\}.$$

\vspace{.1in}

\paragraph{}Just as in the finite type case, one may prove~:

\begin{prop}[\cite{Lus3}] The map $\iota$ is an embedding of $E(\d)$ as a dense open
subset of the (generalized) affine Schubert variety
$$\overline{\mathcal{Y}_{w_0}}:=\{(\mathcal{N}_i)\;|\; \mathcal{N}_i \supset 
\mathcal{L}_{i-1}\quad \text{for\;all}\;i\}.$$
Moreover, each $G_{\d}$-orbit $\mathcal{O}_{\m} \subset E(\d)$ embeds
as a dense open subset of a $P$-orbit $\mathcal{Y}_{w(\m)}$ for some $w(\m) \in
\widehat{\mathfrak{S}}_{\d} \backslash \widehat{\mathfrak{S}}_{r} / 
\widehat{\mathfrak{S}}_{\d}$. 
\end{prop}

\begin{cor} The coefficients of the transition matrix between the canonical
basis $\{\b_{\m}\}$ and the PBW basis $\{\f_{\m}\}$ of $\U^{(1)}_n$
are affine Kazhdan-Lusztig polynomials 
(of type $A^{(1)}_{\sum d_i -1})$.
\end{cor}

\vspace{.2in}

\paragraph{\textbf{Remark.}} Nakajima constructed in \cite{Nak2} 
a similar compactification
of $E(\d)$ for an arbitrary quiver. However, it is only in type A that the coefficients
of the canonical basis are given by Kazhdan-Lusztig polynomials (see, however \cite{BZ2}
for some descriptions of singularities of orbit closures in type $D$ in terms of
\textit{products} of Grassmanians).

\section{Schur-Weyl duality in the geometric setting}

\paragraph{\textbf{7.1.}} Let us first briefly recall the usual Schur-Weyl duality. Let
$V$ be an $n$-dimensional $\C$-vector space. There are two commuting group actions on
$V^{\otimes r}$: $GL(V)$ acts diagonally while $\mathfrak{S}_r$ acts by permutation of
the factors. Moreover we have
$$\C[\mathfrak{S}_r]=\text{End}_{GL(V)}(V^{\otimes r})\qquad \text{if\;}
n \geq r,$$
$$U(\mathfrak{gl}(V)) \twoheadrightarrow \text{End}_{\mathfrak{S}_r}
(V^{\otimes r})$$
The algebra $\text{End}_{\mathfrak{S}_r}
(V^{\otimes r})$ is called the Schur algebra. There are affine and 
quantum versions of this duality. The affine version is obtained by replacing
$\mathfrak{gl}_n$ by $\widehat{\mathfrak{gl}}_n=\mathfrak{gl}_n[z,z^{-1}] 
\oplus \C c$, the symmetric group $\mathfrak{S}_r$ by 
$\widehat{\mathfrak{S}}_r$ and $V$ by $V[z,z^{-1}]$.

\vspace{.1in}

\paragraph{}In the quantum version, we replace $\mathfrak{gl}_n$ by
the quantum group
$\U_v(\mathfrak{gl}_n)$. Recall that this is the $\C(v)$-Hopf algebra generated
by elements $E_i$, $F_i$, $i=1, \ldots, n-1$ and $H_i^{\pm 1}$, $i=1, \ldots ,
n$ satisfying relations
$$H_iH_i^{-1}=H_i^{-1}H_i=1, \qquad H_iH_j=H_jH_i,$$
$$H_iE_j=v^{a_{ij}}E_jH_i, \qquad H_iF_j=v^{-a_{ij}}F_jH_i,$$
$$[E_i,F_j]=
\delta_{ij}\frac{H_iH_{i+1}^{-1}-H_i^{-1}H_{i+1}}{v-v^{-1}},$$
$$[E_i,E_j]=[F_i,F_j]=0\quad \text{if}\; |i-j|>1,$$
$$E_i^{(2)}E_j-E_iE_jE_i+E_jE_i^{(2)}=F_i^{(2)}F_j-F_iF_jF_i+F_jF_i^{(2)}=0\quad
\text{if\;}|i-j|=1.$$
where $(a_{ij})$ is the Cartan matrix of type $A_{n-1}$. The coproduct is
defined as in Section~4.1.

\vspace{.1in}

The standard representation $\C^n$ of $\mathfrak{gl}_n$ admits a quantization $V_q$
(the ``vector'' representation).
Using the coproduct $\Delta$ we define an action of $\U_v(\mathfrak{gl}_n)$
on $V_q^{\otimes r}$. It is well-known that $\U_v(\mathfrak{gl}_n)$ is a
\textit{quasitriangular Hopf algebra}, i.e there exists an element
$\mathcal{R} \in \U_v(\mathfrak{gl}_n) \otimes  \U_v(\mathfrak{gl}_n)$
satisfying the relation $\mathcal{R} \Delta=\Delta^{op}\mathcal{R}$. From
this we construct an action of the Artin braid group on $r$ strands $B_r$
on $V_q^{\otimes r}$ by the formula
$$\sigma_i(v_1 \otimes \cdots \otimes v_r)=v_1 \otimes \cdots \otimes
P \mathcal{R}(v_i \otimes v_{i+1}) \otimes \cdots \otimes v_r,$$
where $\sigma_i$ is the $i$th standard generator of $B_r$ and $P: u \otimes
v \mapsto v \otimes u$ is the permutation (see e.g \cite{ES} for more details on quasitriangular 
Hopf algebras). 
One shows that the action of $B_r$
factorizes through the Hecke algebra $\mathbf{H}_r$, and that, as in the 
classical case, we have

$$\mathbf{H}_r=\text{End}_{\U_v(\mathfrak{gl}_n)}(V_q^{\otimes r})\qquad 
\text{if\;}
n \geq r,$$
$$\U_v(\mathfrak{gl}_n) \twoheadrightarrow \text{End}_{\mathbf{H}_r}
(V_q^{\otimes r})$$
The algebra $\text{End}_{\mathbf{H}_r}
(V_q^{\otimes r})$ is called the $q$-Schur algebra.

\vspace{.1in}

\paragraph{}Finally, the quantum affine version is now obtained by replacing
$\U_v(\mathfrak{gl}_n)$ by $\U_v(\widehat{\mathfrak{gl}}_n)$, the space
$V_q$ by its affine version $V_q[z,z^{-1}]$, and $\mathbf{H}_r$ by
$\widehat{\mathbf{H}}_r$. However, there is a subtle point here: $\U_v
(\widehat{\mathfrak{gl}}_n)$ is the quantization of $\mathfrak{gl}_n[z,z^{-1}]
\oplus \C c$ rather than the quantization of the Kac-Moody algebra $\widehat{\mathfrak{sl}}_n$
or the quantization of the slightly bigger algebra 
$\widehat{\mathfrak{gl}}'_n:=\mathfrak{gl}_n \oplus \bigoplus_{l \in \Z^*}
z^l \mathfrak{sl}_n \oplus \C c$ (unfortunately, this latter algebra is often also denoted
$\widehat{\mathfrak{gl}}_n$ in the litterature). In particular, the Hopf algebra
$\U_v(\widehat{\mathfrak{gl}}_n)$ is best defined through
the \textit{Drinfeld new realization} (see \cite{Dr}) rather than through
quantum Serre-type relations. We refer to \cite{GV} for the precise definition.

\vspace{.2in}

\paragraph{\textbf{Remark.}} In the affine versions of the Schur-Weyl
duality, the modules $V[z,z^{-1}]$ and $V_q[z,z^{-1}]$ are of level zero, i.e
the central element $c$ (or its quantum analog) act trivially. Thus, for the purposes
of Schur-Weyl duality it is
enough to consider the loop algebra $L\mathfrak{gl}_n=
\mathfrak{gl}_n[z,z^{-1}]$ or its quantization.

\vspace{.2in}

\paragraph{\textbf{7.2.}} Following Beilinson-Lusztig-MacPherson in the
finite-type case and Ginzburg-Vasserot and Lusztig in the affine case, we
will now give a geometric construction of the (quantum) Schur-Weyl dualities.
We assume here that $k=\mathbb{F}_q$.

\vspace{.2in}

\paragraph{} Recall the notations of Section~5~: $L$ is a fixed $r$-dimensional
$k$-vector space, $\mathcal{B}$ the flag variety and $\mathcal{F}=
\sqcup_{\d} \mathcal{F}_{\d}$ the $n$-step partial flag variety. For any
$i=1, \ldots, n$ we consider the subvarieties of $\mathcal{F} \times \mathcal{F}$
$$\mathcal{O}^+_i=\{((L_j),(L'_j))\;|L_j=L'_j\;\text{if}\;j \neq i \;
\text{and}\;L'_i \subset L_i,\; {dim}\;L_i/L'_i=1\},$$
$$\mathcal{O}^-_i=\{((L_j),(L'_j))\;|L_j=L'_j\;\text{if}\;j \neq i \;
\text{and}\;L_i \subset L'_i,\; {dim}\;L'_i/L_i=1\},$$
and for any $\d$ we set
$\mathcal{O}^0_{\d}=\Delta \mathcal{F}_{\d} \subset \mathcal{F}_{\d}
\times \mathcal{F}_{\d}$.

\begin{theo}[\cite{BLM}] The assignement $E_i \mapsto \mathbf{1}_{\mathcal{O}_i^{+}}$,
$F_i \mapsto \mathbf{1}_{\mathcal{O}_i^{+}}$, $H_i \mapsto \sum_{\d}
q^{d_i} \mathbf{1}_{\mathcal{O}^0_{\d}}$ extends to a surjective homomorphism
of algebras 
$$\U_v(\mathfrak{gl}_n) \twoheadrightarrow \C_G[\mathcal{F} \times
\mathcal{F}] \otimes_{\C[v,v^{-1}]} \C(v).$$
\end{theo}

\vspace{.2in}

\paragraph{}Next we consider the affine case. Thus $dim\;L=r$, $\mathbb{L}=L((z))$, 
$\mathcal{X}$
is the affine flag variety and $\mathcal{Y}$ is the affine $n$-step partial flag variety. 
Again we set
$$\mathcal{O}^+_i=\{((L_j),(L'_j))\;|L_j=L'_j\;\text{if}\;j \neq i \;
\text{and}\;L'_i \subset L_i,\; {dim}_k\;L_i/L'_i=1\},$$
$$\mathcal{O}^-_i=\{((L_j),(L'_j))\;|L_j=L'_j\;\text{if}\;j \neq i \;
\text{and}\;L_i \subset L'_i,\; {dim}_k\;L'_i/L_i=1\},$$
and
$\mathcal{O}^0_{\d}=\Delta \mathcal{Y}_{\d} \subset \mathcal{Y}_{\d}
\times \mathcal{Y}_{\d}$.

\begin{theo}[\cite{GV}, \cite{Lus3.5}] There exists a surjective algebra homomorphism 
$$\U_v(\widehat{\mathfrak{gl}}_n) \to \C_G[\mathcal{Y} \times \mathcal{Y}] 
\otimes_{\C[v,v^{-1}]}\C(v).$$
 On the subalgebra $U_v(\widehat{\mathfrak{gl}}'_n)$ it is given by
$E_i \mapsto \mathbf{1}_{\mathcal{O}_i^{+}}$,
$F_i \mapsto \mathbf{1}_{\mathcal{O}_i^{+}}$, $H_i \mapsto \sum_{\d}
q^{d_i} \mathbf{1}_{\mathcal{O}^0_{\d}}$.
\end{theo}

\vspace{.2in}

\paragraph{} We are now ready to give the geometric version of Schur-Weyl duality.
The convolution product defines a left action of $\C_G[\mathcal{F}
\times\mathcal{F}]$ on $\C_G[\mathcal{F} \times \mathcal{B}]$ and a right action of
 $\C_G[\mathcal{B}
\times\mathcal{B}]$ on $\C_G[\mathcal{F} \times \mathcal{B}]$. Moreover, it is clear that
these two actions commute. Recall the algebra $\mathbf{S}_{n,r}=\C_G[\mathcal{F} 
\times \mathcal{F}]$ defined in Section~5.
 
\begin{theo}[\cite{GV}] There exists an isomorphism of algebras 
$$\mathbf{S}_{n,r} \simeq 
\mathrm{End}_{\H_r}(V_q^{\otimes r})$$
and an isomorphism of $\mathbf{S}_{n,r}-\H_r$-bimodules
$ \C_G[\mathcal{F} \times \mathcal{B}] \simeq V_q^{\otimes r}$ such that the 
following diagram commutes
\begin{equation}\label{D:1}
\xymatrix{
\mathbf{S}_{n,r}=\C_G[\mathcal{F} \times \mathcal{F}] 
\ar@{~>}[r] \ar[d]_{\sim}&  \C_G[\mathcal{F} \times \mathcal{B}] \ar[d]_{\sim} &
\C_{G}[\mathcal{B} \times \mathcal{B}] \ar@{~>}[l] \ar[d]_{\sim}\\
\mathrm{End}_{\H_r}(V_q^{\otimes r}) \ar@{~>}[r] & V_q^{\otimes r} & \H_r \ar@{~>}[l]
}
\end{equation}
(here $\xymatrix{A \ar@{~>}[r]& B}$ denotes an action of an algebra $A$ on $B$).
\end{theo}

\vspace{.2in}

\paragraph{}In a similar way, we have

\begin{theo}[\cite{VV}] 
There exists an isomorphism of algebras 
$$\widehat{\mathbf{S}}_{n,r} \simeq 
\mathrm{End}_{\AH_r}(V_q[z,z^{-1}])^{\otimes r})$$
and an isomorphism of 
$\widehat{\mathbf{S}}_{n,r}-\AH_r$-bimodules
$ \C_G[\mathcal{Y} \times \mathcal{X}] \simeq (V_q[z,z^{-1}])^{\otimes r}$ 
such that the following diagram
commutes
\begin{equation}\label{D:2}
\xymatrix{
\widehat{\mathbf{S}}_{n,r}=\C_G[\mathcal{Y} \times \mathcal{Y}] 
\ar@{~>}[r] \ar[d]_{\sim}&  \C_G[\mathcal{Y} \times \mathcal{X}] \ar[d]_{\sim} &
\C_{G}[\mathcal{X} \times \mathcal{X}] \ar@{~>}[l] \ar[d]_{\sim}\\
\mathrm{End}_{\AH_r}((V_q[z,z^{-1}])^{\otimes r}) \ar@{~>}[r] & (V_q[z,z^{-1}])^{\otimes r} & 
\AH_r \ar@{~>}[l]
}
\end{equation}
\end{theo}

\vspace{.2in}

\paragraph{\textbf{7.3.}} Using a combination of the constructions in Sections 6 and 7, we can
give a conceptual interpretation of Theorem~4.3. Put
$$\mathcal{O}^-=\{((L_i),(L'_i))\;|L'_i \subset L_i\quad\text{for\;all\;}i\} \subset \mathcal{Y}
\times \mathcal{Y}.$$
Observe that to a point $((L_i),(L'_i)) \in \mathcal{O}^-$ we may associate the representation
of $\Dtn$ on $\bigoplus_{i \in \Z/n\Z} L_i/L'_i$ using the natural maps $\phi_i: L_i/L'_{i}
\to L_{i+1}/L'_{i+1}$ (as in Section~6). Note that such a representation is always nilpotent by
the periodicty condition for affine flags. For a multipartition $\m$ we denote by
$\mathcal{Z}_{\m} \subset \mathcal{O}^-$ the set of all points whose associated representation
of $\Dtn$ is of type $\m$. Let us also denote by $\widehat{\mathbf{S}}_{n,r}^-$ the convolution
algebra $\C_G[\mathcal{O}^-]$.

\begin{prop}[\cite{VV}] The linear map $\varphi: \U_n^{(1)} \to \widehat{\mathbf{S}}_{n,r}^-$
such that $\varphi(\mathbf{1}_{\mathcal{O}_{\m}})= \mathbf{1}_{\mathcal{Z}_{\m}}$ is a surjective
algebra homomorphism.
\end{prop}
\noindent
\textit{Proof.} This follows directly from the definitions of the Hall product and of the
convolution product.\qed

\vspace{.1in}

\paragraph{}To sum up, we have the following diagram of geometric constructions
of algebras :
$$
\xymatrix{
\U_n^{(1)} \ar@{->>}[r] & \C_G[\mathcal{O}^-] \ar@{^{(}->}[r] &
\C_G[\mathcal{Y} \times \mathcal{Y}] \ar@{~>}[r] &
 \C_G[\mathcal{Y} \times \mathcal{X}] &
\C_{G}[\mathcal{X} \times \mathcal{X}] \ar@{~>}[l]
}
$$
On the other hand, the bicommutant Theorem~7.4 gives us a canonical embedding $Z(\AH_r) \subset 
\widehat{\mathbf{S}}_{n,r}$. Recall that by Bernstein's theorem, 
$$Z(\AH_r) = \C[v,v^{-1}]
[X_1^{\pm 1}, \ldots, X_r^{\pm 1}]^{\mathfrak{S}_r}.$$
Let us set $Z^-(\AH_r)=\C[v,v^{-1}]
[X_1^{-1},\cdots ,X_r^{-1}]^{\mathfrak{S}_r}$.

\begin{prop}[\cite{S1}] The morphism $\varphi$ restricts to a surjective algebra homomorphism
$i_r: \mathbf{K} \to Z^-(\AH_r)$.  Moreover, the maps $i_r$ are compatible with the natural
projective system $Z^-(\AH_r) \to Z^-(\AH_{r-1})$ and yield in the projective limit an
isomorphism 
$$i: \mathbf{K} \stackrel{\sim}{\to} Z^-(\AH_{\infty}) \simeq \Gamma \otimes \C[v,v^{-1}].$$
\end{prop}

\paragraph{}In particular, this map gives rise, for each $n$, to two new natural families
of symmetric polynomials, namely $i(\b^*_{(\lambda, \cdots, \lambda)})$ and $i(\pi(\b_{\lambda,
\cdots, \lambda}))$ where $\pi: \U_n^{(1)} \to \mathbf{K}$ is the orthogonal projection. 
When $n=1$ these families are the dual Schur functions and the Schur functions respectively.

\vspace{.1in}

\paragraph{}In the same spirit, we make the following conjecture :

\begin{conj} The restriction of Green's scalar product $\langle\;,\;\rangle_G$ to
$\mathbf{K} \simeq \Gamma \otimes \C[v,v^{-1}]$ is given by the formula
$$(p_\lambda,p_\mu)=\delta_{\lambda\mu}z_{\lambda}q^{(n-1)|\lambda|}(1-q^{-1})^{n|\lambda|}
\prod_{i=1}^{l(\lambda)}\frac{1-q^{n\lambda_i}}{(1-q^{\lambda_i})^2}.$$
\end{conj}

This can be viewed as a generalization of the Hall-Littlewood scalar product.

\vspace{.2in}

\paragraph{\textbf{Remarks.}} i) All the spaces in diagrams (\ref{D:1}) and (\ref{D:2})
are equipped with canonical 
bases. One can show that all the maps are (in some suitable sense) compatible with these
bases (see \cite{VV}, \cite{SV}).\\
ii) Using the map $i: \mathbf{K} \stackrel{\sim}{\to} Z^-(\AH_{\infty})$ and its compatibility 
with canonical bases it is possible
to give some explicit (but very complicated) expressions for the central
elements of $\AH_r$ in terms of the Kazhdan-Lusztig basis (see \cite{S2}). 

\section{Representation theory of Hecke algebras and quantum affine algebras
and graded nilpotent orbits}

\paragraph{}This survey would be very incomplete without mentioning a totally different
aspect of quivers of type $A$ and their canonical bases, namely their interpretations
as Grothendieck rings of finite-dimensional modules for an affine Hecke algebra and for
an affine $q$-Schur algebra.

\vspace{.2in}

\paragraph{\textbf{8.1.}} We will first deal with Hecke algebras. We make a few preliminary 
definitions and remarks. Since $Z(\AH_r)=\C[v,v^{-1}][X_1^{\pm 1}, \ldots, 
X_r^{\pm 1}]^{\mathfrak{S}_r}$,
any character $\theta: Z(\AH_r) \to \C^*$ is given by an $r+1$-tuple $(s_1, \ldots, 
s_r; \epsilon)$ such that $\theta(x_i)=s_i$ and $\theta(v)=\epsilon$. It is convenient to view such
an $r+1$-tuple as a pair $(s,\epsilon) \in T \times \C^*$ where $T \subset GL_r(\C)$ is the
subgroup of diagonal matrices. Note that $\theta$ determines $s$ only up to conjugation
(i.e permutation). By Schur's lemma, $Z(\AH_r)$ acts by a character in any simple 
finite-dimensional module. Let us denote by $Rep_{\theta}(\AH_r)$ the category of 
finite-dimensional $\AH_r$-modules of central character $\theta$.

\vspace{.1in}

\paragraph{}The link between representation
theory of $\AH_r$ and quivers of type $A$ stems from the following fundamental theorem.
 
\begin{theo}[\cite{KL3},\cite{Ginz}] Fix $\theta=(s,\epsilon)$ and assume that $\epsilon$ is not a root of
unity. Set $G_s=Stab_{GL_r(\C)}s$. There is a bijection between the set of simple modules
in $Rep_\theta \AH_r$ and the set of $G_s$-orbits in
$$\mathcal{N}_{\theta}=\{x \in \mathfrak{gl}_r, x\;\text{nilpotent},\; sxs^{-1}=\epsilon^{-2}x
\}.$$
\end{theo}

\paragraph{}The proof of the above theorem (the so-called Deligne-Langlands conjecture for
$\AH_r$) is based on a geometric realization of $\AH_r$ completely different from the one
given in Section~5, as the Grothendieck group of (equivariant) coherent sheaves on a certain
variety $Z \subset T^*\mathcal{B} \times T^*\mathcal{B}$ (see \cite{CG},\cite{KL3}). As a part of this 
realization, two types of modules are constructed for each triple $(s,\epsilon,x)$ : the
simple module $L_{\theta,x}$ and the so-called \textit{standard module} $M_{\theta,x}$,
whose character is easy to compute. The determination of the characters of the simple modules
is now a consequence of the following result.

\begin{theo} Let $\theta=(s,\epsilon)$ and assume that $\epsilon$ is not a root of unity. For
each $x \in \mathcal{N}_{\theta}$ we set $\mathcal{O}_x=G_s\cdot x$. Then for any two
$x,y \in \mathcal{N}_{\theta}$ the multiplicity of $L_{\theta,x}$ in $M_{\theta,y}$ is
$$[M_{\theta,x}:L_{\theta,y}]=\sum_k dim \;\mathcal{H}^k_{\mathcal{O}_y}(IC_{\mathcal{O}_x}).$$
\end{theo}

\vspace{.2in}

\paragraph{}Let us reinterpret the above Theorems in terms of quivers $A_m$. Let us fix
$\epsilon$ (the value of the deformation parameter $v$). For simplicity,
let us denote by 
$Rep_{\epsilon} \AH_r$ the category of all finite-dimensional modules of central charge
$\theta=(s,\epsilon)$ with $Spec(s) \subset \epsilon^{2\Z}$.
For any such $s$, write
$\C^r=\bigoplus_{i \in \Z} V_{\epsilon^{2i}}$ for the corresponding 
eigenspace decomposition of $\C^r$. Now observe that $x(V_{\epsilon^{2i}}) \subset 
V_{\epsilon^{2i-2}}$, i.e $x$ is nothing but a representation on the space $\C^r$
of the quiver $A^\infty_\infty$. Let $\m(s,x)$ be the type of this representation. 
Altogether, Theorems~8.1 and 8.2 are nicely summarized
as follows.

\begin{cor}[\cite{Ar}] Assume that $\epsilon$ is not a root of unity. There exists a linear 
isomorphism
$$\Theta_{\infty}:\bigoplus_{r} Rep^*_\epsilon(\AH_r) \stackrel{\sim}{\to} \U_{\infty,|v=1}$$
such that $\Theta_{\infty}([M_{\theta,x}]^*)= \f_{\m(s,x)}$ and
$\Theta_{\infty}([L_{\theta,x}]^*)=\b_{\m(s,x)}$.\end{cor}

\vspace{.2in}

\paragraph{\textbf{Remark.}} The above Theorem together with the results of Section~6
shows that, in particular, the multiplicities $[M_{\theta,x}:L_{\theta,y}]$ are equal to
evaluation at $1$ of some Kazhdan-Lusztig polynomials of type $A$.

\vspace{.2in}

\paragraph{}The case when $\epsilon^2$ is a primitive $n$th root of unity is more subtle.
In that situation, Theorem~8.2 remains true while in Theorem~8.1 only a certain subset of 
$G_s$-orbits in $\mathcal{N}_{\theta}$ is relevant. Observe that now if 
$Spec(s) \subset \epsilon^{2\Z}$ then to a pair $(s,x)$ is associated a representation of
the cyclic quiver $A_n^{(1)}$ rather than $A^\infty_{\infty}$. The classification of 
simple modules
is best described by the next theorem.

\begin{cor}[\cite{Ar},\cite{Groj}] Assume that $\epsilon^2$ is a primitive $n$th root of 
unity. There exists a linear 
isomorphism
$$\Theta_{n}:\bigoplus_{r} Rep^*_\epsilon(\AH_r)\stackrel{\sim}{\to}\mathbf{C}_{n,|v=1} \subset
\U^{(1)}_{n,|v=1}$$
such that $\Theta_{\infty}([L_{\theta,x}]^*)=\b_{\m(s,x)}$.\end{cor}

\paragraph{}In particular, simple modules in $Rep_\epsilon(\AH_r)$ are indexed by aperiodic 
multipartitions of total dimension $r$. 

\vspace{.1in}

\paragraph{\textbf{Remark.}} i) Ariki proved much more than the results stated above : 
he endowed the space
$\bigoplus_r Rep^*_{\epsilon}(\AH_r)$ with an algebra structure (using the restriction
functor $Rep_{\epsilon}(\AH_{r+r'}) \to Rep_{\epsilon}(\AH_r) \otimes 
Rep_{\epsilon}(\AH_{r'})$) and showed that the map $\Theta$ is actually an algebra 
homomorphism.\\
ii) A consequence of the above Theorem and of the results in Section~6 is that the 
multiplicities $[M_{\theta,x}:L_{\theta,y}]$ are now given by evaluation at $1$ of 
\textit{affine} Kazhdan-Lysztig polynomials of type $A$.\\
iii) Recently, Vasserot succeeded in extending some of the above results to
the case of a \textit{double affine Hecke algebra}, or \textit{Cherednik algebra} 
\cite{Vas2}.

\vspace{.2in}

\paragraph{\textbf{8.2.}} Let us now turn to the case of the quantum affine algebra
$\U_v(\widehat{\mathfrak{sl}}_k)$. It is known that any finite-dimensional representation
of $\U_v(\widehat{\mathfrak{sl}}_k)$ factors through $\widehat{\mathbf{S}}_{k,r}$ for some
$r \gg 0$. Hence, by Schur-Weyl duality we expect a very similar classification result
as in Section~8.1. A geometric construction of $\widehat{\mathbf{S}}_{k,r}$ in the same spirit
as \cite{Ginz}, \cite{KL3} was given in \cite{GV}. 

\begin{theo}[\cite{GV}] Let $\epsilon \in \C^*$. There is a bijection between the simple 
finite-dimensional modules of $\U_\epsilon(\widehat{\mathfrak{gl}}_k)$ and the set of
pairs 
$$\{(s,x)\;|s \in T,\; x \in \mathfrak{gl}_r,\;x^k=0,\; sxs^{-1}=\epsilon^{-2}x\}/GL_r(\C).$$
To each such pair is associated a simple module $L_{s,x}$ and a standard module
$M_{s,x}$, and we have
$$[M_{s,x}:L_{s',y}]=\delta_{s,s'}\sum_l dim\; \mathcal{H}^l_{\mathcal{O}_y}
(IC_{\mathcal{O}_x})$$
where $\mathcal{O}_z$ is the $G_s$-orbit containing $z$.
\end{theo}

\paragraph{}In particular, the classification of simple modules can again be interpreted
in terms of quivers of type $A$.
However, observe that, as opposed to Section~8.1., there are no additional
restrictions when $\epsilon$ is a root of unity. Let $Rep_{\epsilon}
(\U_v(\widehat{\mathfrak{sl}}_k))$ denote the category of all finite-dimensional modules with
a composition series involving only simple modules $L_{s,x}$ with $Spec(s) \subset
 \epsilon^{2\Z}$. For notational convenience, we set $A^{(1)}_\infty=A^\infty_{\infty}$ 
in what follows.

\begin{cor}[\cite{VV}] Let $n \in \N \cup \{\infty\}$ be the order of $\epsilon^2$.
There exists a linear 
isomorphism
$$\Xi_{n}:\bigoplus_{k} Rep^*_\epsilon(\U_v(\widehat{\mathfrak{sl}}_k))
\stackrel{\sim}{\to} \U^{(1)}_{n,|v=1}$$
such that $\Xi_{n}([M_{\theta,x}]^*)= \f_{\m(s,x)}$ and
$\Xi_{n}([L_{\theta,x}]^*)=\b_{\m(s,x)}$.\end{cor}

\paragraph{\textbf{Remark.}} i) Again, the standard modules $M_{s,x}$ can be explicitly
described and their characters computed (see \cite{Vas3}). Also, as for Hecke algebras, the 
multiplicities $[M_{\theta,x}:L_{\theta,y}]$ are evaluation at $1$ of some Kazhdan-Lusztig
polynomials (affine if $\epsilon$ is a root of unity).

\vspace{.2in}

\paragraph{}One may ask for a similar interpretation for Grothendieck groups of categories
of representations of $\U_{\epsilon}(\mathfrak{sl}_k)$ or $\H_r$. Although no geometric
constructions similar to those in \cite{KL3} or \cite{GV} exist for the 
\textit{finite type}
quantum group or Hecke algebra, there exist evaluation maps 
$ev:\U_v(\widehat{\mathfrak{sl}}_k) \to \U_v(\mathfrak{sl}_k)$,
$ev: \AH_r \to \H_r$. These in turn induce maps $ev^*:Rep^*_{\epsilon}
(\U_v(\widehat{\mathfrak{sl}}_k)) \to Rep^*_{\epsilon}(\U_v({\mathfrak{sl}}_k))$
and $ev^* Rep_{\epsilon}^*(\AH_r) \to Rep_{\epsilon}^*(\H_r)$. In fact these maps 
nicely fit in the previous picture, as shown by the following theorems.

\begin{theo}[\cite{Ar},\cite{VV}]
Let $n \in \N \cup \{\infty\}$ be the order of $\epsilon^2$. There is a commutative
diagram
$$
\xymatrix{
\bigoplus_r {Rep}_{\epsilon}^*(\AH_r) 
\ar@{->}[rr]^{\sim}_{\Theta_n} \ar[d]_{ev^*}& & 
\mathbf{C}_{n,v=1} \subset \U^{(1)}_{n|v=1}\ar[d]^\gamma\\
\bigoplus_r Rep_{\epsilon}^*(\H_r) \ar@{->}[rr]^\sim_{\theta_n} &  & V_{\Lambda_0}
}
$$
where $V_{\Lambda_0}$ is the highest-weight integrable representation of $\mathbf{C}_n \simeq
\U^-_v(\widehat{\mathfrak{sl}}_n)$, and $\gamma$ is the action on the highest weight vector.
\end{theo}

\begin{theo}[\cite{VV}] Let $n \in \N \cup \{\infty\}$ be the order of $\epsilon^2$. 
There is a 
commutative diagram
$$
\xymatrix{
\bigoplus_k {Rep}_{\epsilon}^*(\U_v(\widehat{\mathfrak{sl}}_k)) 
\ar@{->}[rr]^{\sim}_{\Xi_n} \ar[d]_{ev^*}& & \U^{(1)}_{n|v=1}\ar[d]^\gamma\\
\bigoplus_k Rep_{\epsilon}^*(\U_v(\mathfrak{sl}_k)) \ar@{->}[rr]^\sim_{\xi_n} &  & 
\Lambda_{\infty|v=1}}
$$
where $\Lambda_{\infty}$ is the Fock space representation of $\U_n$ constructed in
\cite{Hay}, \cite{VV}, and $\gamma$ is the action on the \textit{vacuum vector} of 
$\Lambda_{\infty}$.
\end{theo}

\paragraph{}In addition, the spaces $V_{\Lambda_0}$ and $\Lambda_{\infty}$ are again equipped 
with two bases: a ``canonical basis'' and a ``PBW basis''. These correspond, under the
isomorphisms $\theta_n$ and $\xi_n$ to the (dual) bases of simple modules and standard
modules respectively. Finally, the map $\gamma$ is compatible with these canonical bases.

\vspace{.1in}

\paragraph{}This beautiful interpretation of character formulae in terms of canonical bases
of Hall algebras and Fock spaces or highest weight representations 
was first conjectured by Lascoux-Leclerc-Thibon in \cite{LLT}. It was
proved by Ariki for Hecke algebras (in the more general framework of 
\textit{cyclotomic Hecke algebras}, which
include in particular Hecke algebras of type $B$) and by Varagnolo and Vasserot
for quantum affine algebras. One nice consequence of this construction for quantum affine
algebras is that it gives a new proof of Lusztig's conjecture for the characters of
simple $\U_{\epsilon}(\mathfrak{sl}_k)$-modules when $\epsilon$ is a root of unity.

\section{Some further geometric constructions}

\paragraph{}To conclude this survey, we briefly mention some other interesting
occurences of quivers of type $A$ in geometry.

\vspace{.2in}

\paragraph{\textbf{9.1.Nakajima quiver varieties.}} Nakajima recently introduced a new class
of algebraic varieties and convolution algebras which play an increasingly important role
in geometric representation theory (see [Nak1],[Nak3]). These varieties are defined for an 
arbitrary quiver, but they have a special interpretation
for affine (or finite) Dynkin diagrams as moduli spaces of sheaves on complex surfaces.

\vspace{.1in}

\paragraph{} Let us consider the action of the cyclic group $\Gamma=\Z/n\Z$ on $\C^2$ via the
embedding 
\begin{equation*}
\begin{split}
\Z/n\Z &\to SL(2,\C)\\
1 &\mapsto \begin{pmatrix} \zeta_n &0 \\ 0 &\zeta_n^{-1} \end{pmatrix}
\end{split}
\end{equation*}
where $\zeta_n$ is a primitive $n$th root of unity. This action extends to a natural
action on $\mathbb{P}^2(\C)$, leaving the projective line at infinity $l_{\infty}$
invariant. The McKay correspondence gives a bijection between the set of irreducible 
$\Gamma$-modules and the set of simple roots (vertices) of the Dynkin diagram $A_{n-1}^{(1)}$.
Hence, we may identify the Grothendieck (semi)group $K_0(\mathrm{Rep}\;\Gamma)$ with the
set of dimension
vectors of $\Dtn$. Now let $\mathbf{v}, \mathbf{w}$ be two such dimension vectors, and let
$V$, $W$ be corresponding $\Gamma$-modules. The (smooth) Nakajima 
quiver variety $\mathcal{M}(\mathbf{v},\mathbf{w})$ associated to $\Dtn$ is the moduli space of
pairs $(\mathcal{E}, \phi)$ where $\mathcal{E}$ is a $\Gamma$-equivariant torsion-free sheaf
on $\mathbb{P}^2(\C)$ such that $H^1(\mathbb{P}^2, \mathcal{E}(-l_{\infty})) \simeq V$
and $\phi: \mathcal{E}_{\infty} \stackrel{\sim}{\to} \mathcal{O}_{l_{\infty}} \otimes W$
is a $\Gamma$-invariant isomorphism. Let $\mathcal{M}_{vec}(\mathbf{v},\mathbf{w})
\subset \mathcal{M}(\mathbf{v},\mathbf{w})$ be the subvariety consisiting of pairs
$(\mathcal{E},\phi)$ where $\mathcal{E}$ is a vector bundle. The singular Nakajima variety
is defined as $\mathcal{M}_0(\mathbf{v},\mathbf{w})=\mathcal{M}_{vec}(\mathbf{v},\mathbf{w})
\times (S\C^2)^{\Gamma}$. There is a projective map $\pi:\mathcal{M}(\mathbf{v},\mathbf{w})
\to \mathcal{M}_0(\mathbf{v},\mathbf{w})$ given by $(\mathcal{E},\phi) \mapsto 
(\mathcal{E}^{\star \star},\phi,supp(\mathcal{E}^{\star\star}/\mathcal{E}))$.

\vspace{.1in}

\paragraph{} In the McKay correspondence, the extending vertex of the affine Dynkin diagram
is associated to the trivial representation $\rho_0$. The quiver varieties
$\mathcal{M}(\mathbf{v},\mathbf{w})$ and $\mathcal{M}_0(\mathbf{v},\mathbf{w})$ for
$\Delta_{n-1}$ are defined just as in the affine case with the added condition that $V$ and $W$
do not contain $\rho_0$.

\vspace{.2in}

\paragraph{} Remarkably, the above geometrically defined varieties also admit a very concrete
realization in terms of quivers (for arbitrary type). We give it for type $A_{n-1}^{(1)}$,
the finite type being (as above) obtained from it as a particular case.
Let us consider the \textit{double quiver} of $\Dtn$, where each arrow is replaced
by a pair of arrows going in opposite directions. Let $H$ be the set of
oriented edges of this graph.
For any $h \in H$ we denote by $o(h)$ and $i(h)$ the outgoing
and incoming edges of $h$ respectively.
Let $(\mathbf{v},\mathbf{w})\in \N^n \times \N^n$ with
$\mathbf{v}=(v_a)_{a \in \Z/n\Z}$ and
$\mathbf{w}=(w_a)_{a \in \Z/n\Z}$. Fix some $\Z/n\Z$-graded $\C$-vector spaces
$V=\bigoplus V_a$ and $W=\bigoplus W_a$ such that
$dim\;V=\mathbf{v}$ and ${dim}\;W=\mathbf{w}$. Set
$$E(V,V)=\bigoplus_{h \in H} \mathrm{Hom}(V_{o(h)}, V_{i(h)}),$$
$$L(V,W)= \bigoplus_{k} \mathrm{Hom}(V_k, W_k),\qquad
L(W,V)=\bigoplus_{k} \mathrm{Hom}(W_k, V_k)$$
and
$$\mathbf{M}(\mathbf{v},\mathbf{w})=E(V,V) \oplus L(W,V) \oplus L(V,W).$$
An element of $\mathbf{M}(\mathbf{v}, \mathbf{w})$ will usually be denoted
by its components $(B,i,j)$.

\vspace{.2in}

\paragraph{} Let $\epsilon : H \to \{1,-1\}$ be any function satisfying
$\epsilon(h)+\epsilon(\overline{h})=0$, where $\overline{h}$ is the edge $h$
with opposite orientation. Consider the subvariety
$$\mathbf{M}^0(\mathbf{v},\mathbf{w})=\{(B,i,j) \in \mathbf{M}(\mathbf{v},\mathbf{w})\;|
 \sum_{o(h)=k} \epsilon(h)B_{\overline{h}}B_h+i_kj_k =0\;\text{for\;all}\;k\}.$$
The group $G_V=\prod_k GL(V_k)$ acts on $\mathbf{M}^0(\mathbf{v},
\mathbf{w})$ by conjugation.
The quiver varieties are obtained as quotients of the variety by $G_V$. According to the
Geometric Invariant Theory philosophy, there are essentially two ways of doing this.
Define an element
$(B,i,j)$ to be \textit{semistable} if the following condition holds~: the only $B$-invariant,
graded subspace $S \subset \bigoplus V_k$ contained in $\mathrm{Ker}\;j$ is $\{0\}$. Let
$\mathcal{M}^{0,ss}(\mathbf{v},\mathbf{w}) \subset
\mathbf{M}^0(\mathbf{v},\mathbf{w})$ denote the open subvariety of semistable points.

\begin{theo}[\cite{Nak0}, \cite{VV2}] We have $\mathcal{M}_0(\mathbf{v},\mathbf{w})=
\mathbf{M}^0(\mathbf{v},\mathbf{w})//G_V$ (categorical quotient) and 
$\mathcal{M}(\mathbf{v},\mathbf{w})=
\mathbf{M}^{0,ss}(\mathbf{v},\mathbf{w})/G_V$ (geometric quotient).
\end{theo}

\vspace{.2in}

\paragraph{} Quite remarkably, in type $A$ the quiver varieties are again related to partial
flag varieties. Let us consider the case of the quiver $\Delta_n$.
Let $\mathbf{v},\mathbf{w}$ be $n$-tuples of positive integers, and define another $n$-tuple
$\mathbf{a}$ as follows~:
$$a_1=w_1 + \cdots + w_{n-1} -v_1,\qquad a_n=v_{n-1},$$
$$a_i=w_i + \cdots + w_{n-1} -v_i+v_{i-1},\qquad \text{for\;}i=2, \ldots, n-1.$$
Behold the partial
flag variety
$$\mathcal{F}_{\mathbf{a}}=\{U_{1} \subset U_2 \subset \cdots \subset U_{n-1} 
\subset \C^r\;|\;dim\;U_i/U_{i-1}=a_{i}\},$$
where $r=\sum_i a_i$.
The cotangent bundle 
$T^*(\mathcal{F}_{\mathbf{a}})$ can be identified with the variety of pairs
$$\{((U_i),x)\;|x(U_i) \subset U_{i+1}\} \subset \mathcal{F}_{\mathbf{a}}
\times \mathcal{N},$$
where $\mathcal{N} \subset \mathfrak{gl}_r$ denotes the set of nilpotent elements. It is
known that the image of the natural projection $p:T^*\mathcal{F}_{\mathbf{a}} \to
\mathcal{N}$ is the Zariski closure of one orbit, say $\mathcal{O}_{\mathbf{a}}$.

\vspace{.1in}

\paragraph{}Now let $y \in \mathcal{N}$ be a nilpotent element with $w_i$ Jordan blocks of 
length $i$ and let $z,h \in \mathfrak{gl}_r$
be such that $(y,h,z)$ is an $\mathfrak{sl}_2$-triple. The \textit{Slodowy slice}
at $y$ is defined as
$$\mathcal{S}_y =\{ u \in \mathcal{N}\;| [u-y,z]=0\}.$$
The variety $\mathcal{S}_y$ is transverse to all nilpotent orbits containing $y$ in their 
closure, and plays an important role in the Brieskorn-Slodowy approach to Kleinian
singularities. The following was conjectured by Nakajima and proved by Maffei 
\cite{Maffei}.

\begin{theo} We have $$\mathcal{M}_0(\mathbf{v},\mathbf{w})\simeq \mathcal{S}_y \cap
\overline{\mathcal{O}_{\mathbf{a}}}$$
 and 
$$\mathcal{M}(\mathbf{v},\mathbf{w}) \simeq
p^{-1}(\mathcal{S}_y \cap
\overline{\mathcal{O}_{\mathbf{a}}}) \subset T^*\mathcal{F}_{\mathbf{a}}.$$
\end{theo}

\vspace{.2in}

\paragraph{}We note that a similar link between quiver varieties and 
transversal slices to
nilpotent orbits appears in \cite{MV}, and is related there to some infinite
Grassmannians.

\vspace{.2in}

\paragraph{\textbf{9.2.Parabolic structures.}} Let $k$ be an algebraically closed field
and let $X$ be a curve defined over $k$. If $x \in X$ is a smooth point then the category
of torsion sheaves (or skyscraper sheaves) on $X$ supported at $x$ is equivalent to the
category of finite-dimensional representations of the local ring $\mathcal{O}_x$ at $x$.
Since this is a discrete valuation ring, this category is equivalent to 
$Mod_k(\Delta^{(1)}_1)$. In particular, the Hall algebra of the category of
coherent sheaves on $X$ (over $\mathbb{F}_q$) contains a copy of $\U^{(1)}_1 \simeq \Gamma
\otimes \C[v,v^{-1}]$ at each (smooth) point.

\vspace{.1in}

\paragraph{}Now let us assume in addition that $X$ is endowed with an action of
the cyclic group $G=\Z/n\Z$, and that $x$ is an isolated fixed point for this action.
Hence $G$ acts on $\mathcal{O}_x$ as $i\cdot f=\zeta_n^i f$, where $\zeta_n=e^{2i\pi/n}$.
In that situation, the category of $G$-equivariant torsion 
sheaves supported at $x$ is equivalent to $Mod_k(\Dtn)$. 

\vspace{.1in}

\paragraph{}Finally, let us consider the case of a finite group $G$ acting on $X$, and let
$\pi: X \to X/G \simeq Y$ denote the quotient map. Assume that this map is ramified at points
$y_1, \ldots, y_r$ with ramification indices $p_1, \ldots, p_r$. Then, locally over each
$y_i$ we are in the situation just described with $G=\Z/d_i\Z$. In particular, the category
of $G$-equivariant coherent sheaves on $X$ contains a copy of $Mod_k(\Delta^{(1)}_{d_i})$
for each $i$. This has applications in [S3] to the study of Hall algebras of such categories (see also \cite{GL}).

\vspace{.2in}

\paragraph{}Closely related to the above is the notion of \textit{parabolic structure}.
Let $X$ be a curve and let $x \in X$ be a smooth point. A $p$-parabolic structure at $x$
on a coherent sheaf $\mathcal{F}$ on $X$ is a flag of coherent sheaves
$$\mathcal{F}=\mathcal{F}_0 \subset \cdots \subset \mathcal{F}_p=\mathcal{F}(-x).$$
Note that if $\mathcal{F}$ is locally free at $x$ then we have an induced filtration
$$(\mathcal{F}/\mathcal{F}_0)_{|x} 
\subset \cdots \subset (\mathcal{F}_p/\mathcal{F}_0)_{|x}
\simeq \mathcal{F}_{|x},$$ i.e a parabolic flag in the fiber of $\mathcal{F}$ at $x$.
This notion can be seen as some globalization of the concept of an affine flag as
in Section~5. In \cite{FK} the authors consider a similar situation 
in the case of a surface : using a variant of the Hall algebra, they construct an action
of $\widehat{\mathfrak{gl}}_p$ on the cohomology of certain moduli spaces of coherent
sheaves on surfaces with $p$-parabolic structure along a rational curve.

\small{

\vspace{4mm}
Olivier Schiffmann,\texttt{schiffma@dma.ens.fr},\\
DMA, \'Ecole Normale Sup\'erieure, 45 rue d'Ulm, 75230 Paris Cedex 05-FRANCE
\end{document}